\documentclass{article}
\usepackage[T2A]{fontenc}
\usepackage[cp1251]{inputenc}
\usepackage[english,russian]{babel}
\usepackage[tbtags]{amsmath}
\usepackage{amsfonts,amssymb,mathrsfs,amscd}
\usepackage{hyperref}
\numberwithin{equation}{section}
\newtheorem{theorem}{Теорема}
\newtheorem{lemma}{Лемма}[section]

\def\const{\mathrm{const}}
\def\Expect{\mathsf{E}}
\def\cov{\mathrm{cov}}
\def\Ext{\mathrm{Ext}}
\def\dist{\mathrm{dist}}
\makeatletter
\def\blfootnote{\xdef\@thefnmark{}\@footnotetext}
\makeatother
\begin{document}

\title{Малые подграфы и их расширения в случайном дистанционном
графе}
\author{А.\,В.~Буркин, М.\,Е.~Жуковский}
\date{}
\maketitle

\begin{abstract}
  В настоящей работе доказываются утверждения, касающиеся распределения малых
  подграфов в последовательности случайных дистанционных графов. Ранее было доказано утверждение
  о пороговой вероятности для свойства содержать фиксированный строго
  сбалансированный граф, в этой же статье мы получаем
  более сильные обобщения этого результата.
\end{abstract}

\markright{Малые подграфы и расширения в случайном дистанционном графе}

\blfootnote{Настоящая работа выполнена при финансовой поддержке Министерства образования и науки РФ по Программе повышения конкурентоспособности РУДН среди ведущих мировых научно-образовательных центров на 2016--2020 гг., а также грантов РФФИ №~15-01-03530, РФФИ №~16-31-60052.}

\allowdisplaybreaks[3]

\section{Введение и история задачи}
\label{intro}

П.~Эрдешем и А.~Реньи в 1959--1960 гг. была предложена модель
случайного графа $G(n, p)$, в которой каждое ребро присутствует с
вероятностью $p$ независимо от остальных ребер (см.
\cite{erdos_renyi_1959}, \cite{erdos_renyi_1960}). Иными словами, $G(n,
p)$ есть случайный элемент со значениями в множестве $\Omega_n$ всех
неориентированных графов $G = (V_n, E)$
с множеством вершин $V_n = \{1, \ldots, n \} $ без петель и кратных ребер и распределением на $ \mathcal{F}_n =
2^{\Omega_n} $, заданным формулой $ \Prob(G) = p^{|E|} (1 - p)^{C_n^2 -
|E|}$.

В основополагающих работах П.~Эрдешем и А.~Реньи был поставлен
вопрос о распределении малых подграфов в случайном графе $G(n,p)$.
Позже этой задачей занимались Б.~Боллобаш~\cite{bollobas_1981}, А.~Ручински, Э.~Винс~\cite{rucinski_vince}, Дж.~Спенсер~\cite{spencer} и др.
Монографии~\cite{bollobas}--\cite{probmethod} посвящены более
полному обзору результатов о распределении малых подграфов в
случайном графе Эрдеша--Реньи и описанию других его
асимптотических свойств. В настоящей работе мы получили ряд
результатов об асимптотическом распределении малых подграфов в другой модели случайного графа, называемой случайным
дистанционным графом, определение которой будет дано в следующем
разделе. Задачами такого типа занимались А.\,Р.~Ярмухаметов (см.,
например, \cite{yarm_conn}), М.\,Е.~Жуковский (см., например,
\cite{zhukovskii_copies}, \cite{Zhuk_ppi}), С.\,Н.~Попова
(см.~\cite{Sveta}).

Далее в этом разделе мы сформулируем некоторые результаты,
относящиеся к асимптотическим свойствам случайного графа
Эрдеша--Реньи.

Пусть $\mathcal{A} = \mathcal{A}(n)$ --- произвольное свойство графов.
{Пороговой вероятностью} свойства $\mathcal{A}$ для случайного графа $G(n,p)$
называется такая функция $p^* = p^*(n)$, что
$\lim_{n\rightarrow\infty}\Prob(\mathcal{A})=0$ при $p =
o(p^*)$, $n \to \infty$, и $\lim_{n\rightarrow\infty}\Prob(\mathcal{A})=1$ при $p = w(p^*)$, $n \to \infty$ (или
наоборот). Здесь $f(n) = o(g(n))$ ($f(n) = w(g(n))$) означает, что
для любого $C > 0$ существует $n_0 > 0$, такое, что для любого $n
> n_0$ выполнено $|f(n)| < C |g(n)|$ ($C|g(n)| < |f(n)|$). Для
этих отношений мы также будем использовать обозначения $f(n) \ll
g(n)$ и $f(n) \gg g(n)$ соответственно. Функция $p^* = p^*(n)$
называется \emph{точной пороговой вероятностью}, если
$\lim_{n\rightarrow\infty}\Prob(\mathcal{A})= 0$ при $p \le
c p^*$ для некоторого $c < 1$, и
$\lim_{n\rightarrow\infty}\Prob(\mathcal{A})= 1$ при $p \ge
c p^*$ для некоторого $c > 1$ (или наоборот).

Для произвольного графа $F$ будем обозначать $v(F)$ и $e(F)$
количество вершин и количество ребер соответственно. Напомним, что
граф $F$ называется \emph{строго сбалансированным}, если
\[
    \rho(H) < \rho(F)
\]
для любого собственного непустого подграфа $H \subset F$, где
$\rho(H) = e(H)/v(H)$ --- \emph{плотность} графа $H$. Граф
называется \emph{сбалансированным}, если неравенство нестрогое.

\emph{Максимальной плотностью} графа $F$ называется величина
\[
    \rho^{\max}(F) = \max_{\substack{H \subseteq F \\
    v(H) \ne 0}} \frac{e(H)}{v(H)}.
\]

П.~Эрдеш и А.~Реньи доказали в \cite{erdos_renyi_1960} теорему о
пороговой вероятности для свойства содержать связный
сбалансированный граф. Это утверждение было в 1981 году обобщено
Б.~Боллобашем на случай произвольного графа
(см.~\cite{bollobas_1981}, в 1985 году А.~Ручински и Э.~Винс
опубликовали  более простое доказательство \cite{rucinski_vince}).

\begin{theorem}[Б.~Боллобаш; А.~Ручински, Э.~Винс]
\label{copies_thrld_erdos}
Пусть $F$ --- произвольный фиксированный граф.
Тогда функция $p^*=n^{-1/\rho^{\max}(F)}$ является пороговой
вероятностью свойства содержать копию $F$ для случайного графа
$G(n,p)$. Выполнен также закон больших чисел для числа копий $X_F$
графа $F$ в $G(n,p)$: при $p\gg p^*$ для любого $\varepsilon
> 0$
\[
    \Prob\left(\left|\frac{X_F}{\Expect X_F} - 1 \right|
    < \varepsilon \right) \to 1.
\]
\end{theorem}

Б.~Боллобашем в \cite{bollobas_1981} было найдено асимптотическое
распределение числа копий $X_F$ строго сбалансированного графа $F$
в $G(n, p)$, если $p$ --- пороговая вероятность свойства содержать
граф $F$.

\begin{theorem}[Б.~Боллобаш] \label{poisson_classic}
Пусть $F$ --- строго сбалансированный граф с $k$ вершинами и $l$
ребрами и $a$ --- количество его автоморфизмов. Пусть $p \sim c
n^{-k/l}$, $c > 0$. Тогда распределение величины $X_F$ слабо сходится к
пуассоновскому с параметром $\lambda = c^l / a$.
\end{theorem}

Перейдем, наконец, к описанию результата, полученного Дж.~Спенсером. Речь пойдет о так называемых свойствах расширений.

Пусть $H$ --- граф с вершинами $z_1, \ldots, z_d, y_1, \ldots,
y_k$, где $R = \{z_1, \ldots, z_d\}$ --- множество \emph{корней}.
\emph{Сетью} называется пара $(R, H)$. Говорят, что граф $G$
удовлетворяет \emph{свойству расширения} $\Ext(R, H)$, если для
любых $v_1, \ldots, v_d \in V(G)$ найдутся такие $w_1, \ldots, w_k
\in V(G)$, что $\{z_i, y_j\} \in E(H) \Rightarrow \{v_i, w_j\} \in
E(G)$ для любых $i\in\{1,\ldots,d\}$, $j\in\{1,\ldots,k\}$ и
$\{y_i, y_j\} \in E(H) \Rightarrow \{w_i, w_j\} \in E(G)$ для
любых $i,j\in\{1,\ldots,k\}$.

Пусть $l = e(H) - e\left(H|_R\right)$, где $H|_R$ --- подграф $H$,
индуцированный на множестве $R$. В общем случае величины $k$ и $l$
будем обозначать $v(R, H)$ и $e(R, H)$ соответственно. Величина
$\rho(R, H) = l/k$ называется \emph{плотностью} сети $(R, H)$.
\emph{Подсетью} называется сеть $(R, S) = \left(R, H|_S\right)$,
где $R \subset S \subseteq
V(H)$. В \emph{собственной} подсети $S \ne V(H)$. Сеть $(R, H)$
называется \emph{строго сбалансированной}, если $\rho(R, S) <
\rho(R, H)$ для всех собственных подсетей $(R, S)$. Она
называется \emph{сбалансированной}, если неравенства нестрогие.
Сеть $(R, H)$ называют \emph{нетривиальной}, если каждая корневая
вершина $z$ соединена ребром в $H$ с хотя бы одной вершиной $y \in
V(H) \setminus R$.

Дж.~Спенсером в \cite{spencer} была доказана следующая теорема.

\begin{theorem}[Дж.~Спенсер]
\label{extension_erdos} Пусть $(R, H)$ --- нетривиальная строго
сбалансированная сеть. Тогда существуют такие числа
$0<\varepsilon<K$, что
$$
 \text{если } p\le\varepsilon n^{-k/l}(\ln n)^{1/l}, \text{ то }
 \lim_{n\rightarrow\infty}\Prob(\Ext(R, H))=0;
$$
$$
 \text{если } p\ge K n^{-k/l}(\ln n)^{1/l}, \text{ то }
 \lim_{n\rightarrow\infty}\Prob(\Ext(R, H))=1.
$$
Пусть $c_1$ есть число автоморфизмов графа $H$, оставляющих корни
на своих местах. Пусть, кроме того, $c_2$ обозначает количество
биективных отображений $R$ на себя, которые можно продолжить до
некоторого автоморфизма $H$. Если $\lambda = \const > 0$ и для $p
= p(n)$ выполнено
\[
    n^k p^l / c_1 = \ln\left(n^d / (c_2 \lambda)\right),
\]
то
\[
    \lim_{n\rightarrow\infty}\Prob(\Ext(R, H))= e^{-\lambda}.
\]
\end{theorem}

В \cite{spencer} доказано также обобщение первой части данной
теоремы (существование пороговой вероятности) на случай
произвольной сети $(R,H)$.

В следующем разделе мы определим \emph{случайный дистанционный
граф} и приведем формулировки доказанных нами теорем для этой
модели, аналогичных теоремам \ref{copies_thrld_erdos}--\ref{extension_erdos}.

\section{Описание модели и новые результаты}

В настоящей работе рассматривается \emph{(симметричный) полный
дистанционный граф}
\[
    G = G(n, n/2, n/4) = (V, E),\quad
    (n \equiv 0\mod 4),
\]
в котором
\[
    V = \{{\bf x} = (x_1, \dots, x_n) \colon x_i \in \{ 0, 1 \},\ x_1 + \ldots + x_n = n/2 \},
\]
\[
    E = \{ \{ {\bf x}, {\bf y} \} \colon \langle{\bf x}, {\bf y}\rangle = n/4 \},
\]
где $\langle{\bf x}, {\bf y}\rangle$ обозначает евклидово
скалярное произведение.

Этот граф называется дистанционным, поскольку его ребра
соответствуют парам вершин, находящихся на определенном расстоянии
друг от друга. Рассмотрение дистанционных графов мотивировано классической
задачей комбинаторной геометрии о хроматическом числе пространства (см.~\cite{raigor_linalg} и \cite{raigor_borsuk}). Впервые дистанционный граф $G(n,r,s)$ (в нашем случае $r=n/2$, $s=n/4$) рассмотрели в 1981 году П.~Франкл и Р.\,М.~Уилсон. С помощью этого графа они показали, что хроматическое
число пространства $\mathbb{R}^n$ растет экспоненциально (см.~\cite{FW}).
В 1991 году Дж.~Кан и Г.~Калаи применили результаты Франкла и Уилсона для
опровержения классической гипотезы Борсука (см.~\cite{raigor_linalg} и \cite{KK}).
Таким образом, изучение внутренней структуры дистанционного графа и его
подграфов играет исключительно важную роль. Сейчас с исследованием
дистанционных графов связаны одни из самых широко изучаемых разделов
комбинаторной геометрии (см.~\cite{raigor_linalg}, \cite{raigor_borsuk},
\cite{Moser}).

Количество вершин этого графа будем обозначать $N = N(n)$, а его
степень (граф, очевидно, является регулярным) --- $N_1 = N_1(n)$.
Заметим, что в силу формулы Стирлинга
\[
    N = C_n^{n/2} \sim \sqrt{\frac{2}{\pi}} \cdot \frac{2^n}{\sqrt{n}},\
    N_1 = \left(C_{n/2}^{n/4}\right)^2 \sim \frac{4}{\pi}\cdot \frac{2^n}{n}.
\]

Нас интересует \emph{случайный дистанционный граф} $G_p = G_p(n,
n/2, n/4)$ --- случайный подграф $G$, в котором каждое ребро
полного дистанционного графа содержится с вероятностью $p = p(n)$
независимо от других ребер \big(т.е. $G_p$ --- случайный элемент со значениями в множестве $\Omega_n^{\dist}$ всех остовных подграфов $G' = (V, E')$ графа $G$ и распределением на $\mathcal{F}_n^{\dist} = 2^{\Omega_n^{\dist}}$, заданным формулой $\Prob(G') = p^{|E'|} (1 - p)^{|E| - |E'|}$\big). Случайные подграфы широко применяются в
вероятностном методе (см., например, \cite{probmethod}). Асимптотические
свойства случайного дистанционного графа изучались, например, в
работах~\cite{yarm_conn}--\cite{Sveta},~\cite{uspekhi}. В
\cite{zhukovskii_copies} была доказана следующая теорема, являющаяся аналогом
теоремы П.~Эрдеша и А.~Реньи о пороговой вероятности для свойства
содержать копию строго сбалансированного графа (которая, в свою
очередь, является частным случаем
теоремы~\ref{copies_thrld_erdos}).

\begin{theorem}[М.\,Е.~Жуковский] \label{zhukovskii_thrld}
Пусть $F$ --- строго сбалансированный граф с $k$ вершинами и $l$
ребрами. Тогда функция
\[
    p^* = N^{-k/l} \sqrt{\ln N}
\]
является пороговой вероятностью свойства содержать копию графа $F$
для случайного графа $G_p$.
\end{theorem}

\subsection{Новые результаты}
\label{new_results}

В настоящей статье мы обобщаем теорему \ref{zhukovskii_thrld} на
случай произвольного графа. Здесь мы используем обозначение $X_F$
для числа копий графа $F$ в случайном графе $G_p$. Заметим, что пороговые
вероятности в случае произвольного случайного подграфа определяются так же,
как и в случае $G(n, p)$.

\begin{theorem} \label{copies_thrld}
Пусть $F$ --- произвольный фиксированный граф.
Тогда функция
\[
    p^*=N^{-1/\rho^{\max}(F)} \sqrt{\ln N}
\]
является пороговой вероятностью свойства содержать копию $F$ для
случайного графа $G_p$. При $p\gg p^*$ для любого $\varepsilon
> 0$
\[
    \Prob\left(\left|\frac{X_F}{\Expect X_F} - 1 \right| <
    \varepsilon \right) \to 1.
\]
\end{theorem}

Теорема будет доказана в разделе \ref{copies_thrld_proof}.

Заметим, что в силу теоремы~\ref{copies_thrld} пороговая
вероятность $p^*$ имеет асимптотику
\[
 p^*\asymp N^{-1/\rho^{\max}(F)}\frac{N}{N_1}
\]
(для двух стремящихся к бесконечности последовательностей $f(n)$ и
$g(n)$ мы пишем $f(n)\asymp g(n)$, если существуют такие числа
$0<c<C$, что для любого $n\in\mathbb{N}$ выполнено
$c|g(n)|\le|f(n)|\le C |g(n)|$) и, тем самым, совпадает с
пороговой вероятностью из теоремы~\ref{copies_thrld_erdos} (в
полном графе $K_n$ степень любой вершины равна $n-1$).

Мы также доказали теорему, аналогичную теореме
\ref{poisson_classic}.

\begin{theorem} \label{poisson}
Пусть $F$ --- строго сбалансированный граф с $k$ вершинами и $l$
ребрами и $a$ есть число автоморфизмов $F$. Пусть
\[
    p \sim c N^{-k/l} \frac{N}{N_1},
\]
где $c = \const > 0$. Тогда распределение $X_F$ слабо сходится к
пуассоновскому с параметром $\lambda = c^l / a$.
\end{theorem}

Она будет доказана в разделе \ref{poisson_proof}.\\

Обратимся, наконец, к свойствам расширений. Такие свойства
являются монотонными (см., например,~\cite{janson_rg}) и поэтому
для них существуют пороговые вероятности (см.~\cite{janson_rg}).
Тем не менее, для многих сетей $(R,H)$ и для любых $p$ свойства
$\Ext(R, H)$ с вероятностями, стремящимися к 1, не выполнены для
некоторых подпоследовательностей случайных дистанционных графов в
силу разреженности дистанционного графа $G(n,n/2,n/4)$. В
частности, если $n$ не делится на 8, то (см.,
например,~\cite{Zhuk_matsbornik}) в графе $G(n,n/2,n/4)$ найдутся
три вершины, не обладающие общим соседом (в данном случае
рассматривается следующее свойство расширения: любые три вершины
обладают общим соседом). В то же время при $8|n$ в этом графе у
любых трех вершин найдется достаточно большое количество соседей.
Поэтому для подобных свойств расширений пороговую вероятность не
удается представить в удобном виде, как
это сделано в теореме~\ref{extension_erdos} для случайного графа
$G(n,p)$. Такая проблема возникает, очевидно, из-за того, что
в качестве множества корней можно взять любой набор $d$ вершин из $V$,
а среди таких наборов встречаются комбинации, приводящие к
``исключениям''. Естественным решением является сузить систему множеств
корней. Оказывается, это можно сделать так, чтобы мощность получившейся
системы была асимптотически равна мощности семейства всех наборов вершин.
Таким образом, при этих ограничениях мы не теряем много информации, и
полученные новые свойства расширений достаточно аккуратно отражают
структуру графа. Итак, определим эти свойства.

Пусть $f(n)$ --- произвольная последовательность положительных
чисел. Пусть, кроме того, $(R, H)$ --- нетривиальная строго
сбалансированная сеть с $V(H)=\{z_1, \ldots, z_d, y_1,\ldots,
y_k\}$, $R = \{z_1, \ldots, z_d\}$ и $e(R, H) = l$. Рассмотрим
произвольные вершины $\mathbf{v}^1=(v_1^1,\ldots,v_n^1), \ldots,
\mathbf{v}^d=(v_1^d,\ldots,v_n^d) \in V$. Напомним, что вершины нашего
графа находятся в пространстве $\{0, 1\}^n$. Обозначим
$\delta_1,\ldots,\delta_{2^d}\in\{0,1\}^d$ различные
$d$-последовательности из нулей и единиц, упорядоченные
лексикографически:
$\delta_1=(1,\ldots,1)>\ldots>(0,\ldots,0)=\delta_{2^d}$. Разобьем
множество $\{1, \ldots, n\}$ на подмножества $B_1, \ldots,
B_{2^d}$ следующим образом: $i\in B_j$ тогда и только тогда, когда
$(v_i^1,\ldots,v_i^d)=\delta_j$.
Положим
\[
 x_j=x_j\left(\mathbf{v}^1,\ldots,\mathbf{v}^d\right)=
 |B_j|-\left[n/2^d\right] \text{ при }
 j\in\left\{1,\ldots,2^{d}-1\right\},
\]
\begin{equation}
 x_{2^d}=x_{2^d}\left(\mathbf{v}^1,\ldots,\mathbf{v}^d\right)=
 |B_{2^d}|-n+\left(2^{d}-1\right)\left[n/2^d\right],
\label{x}
\end{equation}
где $[\cdot]$ --- целая часть числа.
Обозначим $\tilde V_f^d$ множество всех $d$-последователь\-но\-стей
вершин из $V$, для которых $|x_j| \le f(n)$, $j \in \{1, \ldots,
2^d\}$. Будем говорить, что остовный подграф $G'$ дистанционного графа
$G$ обладает свойством $\Ext_f^{\dist}(R, H)$, если для любых
$(\mathbf{v}^1, \ldots, \mathbf{v}^d) \in \tilde V_f^d$ найдутся
такие $\mathbf{w}^1,$ $\ldots, \mathbf{w}^k \in V$, что $\{z_i,
y_j\} \in E(H) \Rightarrow \{\mathbf{v}^i, \mathbf{w}^j\} \in
E(G')$ для любых $i\in\{1,\ldots,d\}$, $j\in\{1,\ldots,k\}$ и
$\{y_i, y_j\} \in E(H) \Rightarrow \{\mathbf{w}^i, \mathbf{w}^j\}
\in E(G')$ для любых $i,j\in\{1,\ldots,k\}$. Иными словами,
свойство $\Ext_f^{\dist}(R, H)$ получается из $\Ext(R, H)$
рассмотрением лишь тех $d$-последовательностей вершин из $V$,
которые принадлежат $\tilde V_f^d$. Таким образом, нас
интересуют только последовательности вершин, разбивающие $\{1,
\ldots, n\}$ на приблизительно равные подмножества. Несколько позже
мы увидим, что $|\tilde V_f^d| \sim |V|^d$. Теорема, сформулированная
ниже, выполнена при условии $f\ll n^{2/3}$ (на самом деле условие на
$f$ можно ослабить, но для наших целей это не принципиально).

\begin{theorem} \label{extdist}
Пусть $c_1$ есть число автоморфизмов $H$, которые оставляют
на месте каждый корень $z_i \in R$. Пусть $p = p(n)$ удовлетворяет равенству
\[
    N^k \left(\frac{N_1}{N}\right)^l p^l / c_1 = d \ln N.
\]
Тогда $p$ является для свойства
$\Ext_f^{\dist}(R, H)$ точной пороговой вероятностью.
\end{theorem}

Теорема \ref{extdist} будет доказана в разделе
\ref{extdist_proof}. Перед доказательством этих теорем мы в
разделе~\ref{lemmas_proofs} докажем вспомогательные леммы,
сформулированные в разделе~\ref{lemmas_statement}, которые
представляют и самостоятельный интерес.

\subsection{Вспомогательные утверждения}
\label{lemmas_statement}

В этом разделе мы сформулируем несколько утверждений, касающиеся
полного дистанционного графа $G$.

Пусть $f = f(n)$ --- произвольная последовательность положительных
чисел, причем $f\ll n^{2/3}$. Пусть, кроме того, $(R, H)$ ---
произвольная сеть с $V(H)=\{z_1, \ldots, z_d, y_1,\ldots, y_k\}$,
$R = \{z_1, \ldots, z_d\}$ и $e(R, H) = l$. Для $\mathbf{v}^1,
\ldots,$ $\mathbf{v}^d \in V$ обозначим $M_{(R, H)} (\mathbf{v}^1,
\ldots, \mathbf{v}^d)$ количество инъективных отображений из $V(H)$ в
$V$, переводящих $z_i$ в $\mathbf{v}^i$, $i\in\{1,\ldots,d\}$, и
сохраняющих ребра между вершинами, среди которых хотя бы одна не
является корнем. Поскольку величина $M_{(R, H)} (\mathbf{v}^1,
\ldots, \mathbf{v}^d)$ не зависит от выбора конкретных вершин, а лишь от значений
$|B_{j}|$ (см. раздел~\ref{new_results}), $j \in \{1, \ldots,
2^d\}$, которые, в свою очередь, задаются числами
$x_1,\ldots,x_{2^d}$, будем обозначать $M_{(R,
H)}^{\vec{x}}=M_{(R, H)} (\mathbf{v}^1, \ldots, \mathbf{v}^d)$,
где вектор $\vec{x} = (x_1, \ldots, x_{2^d})$ определен
в~(\ref{x}).

\begin{lemma} \label{ext1}
Найдется такая функция $M_{(R, H)}'=M_{(R, H)}'(n)$, не зависящая от $\vec{x}$, что $M_{(R,
H)}^{\vec{x}} = M_{(R, H)}' \left(1 + O\left(f(n) n^{-0.8}\right)\right)$ при $n\rightarrow\infty$ равномерно по всем $\vec{x}$ с условием $|x_j| \le f(n)$, $j\in\{1,\ldots,2^d\}$.
\end{lemma}

Лемма \ref{ext1} будет доказана в разделе~\ref{ext1_proof}.

В следующем утверждении мы получили асимптотику количества
вхождений произвольного графа $F$ в дистанционный граф $G$.

\begin{lemma} \label{copies_lem}
Пусть $M_F$ --- количество мономорфизмов графа $F$ с $k$ вершинами
и $l$ ребрами в $G$. Тогда
\[
    M_F \sim M(k, l) := N^k \left(\frac{N_1}{N}\right)^l.
\]
\end{lemma}

В то время как коротким доказательство леммы \ref{ext1} не назовешь, мы нашли
элегантное и лаконичное доказательство леммы \ref{copies_lem}, которое
использует индукцию по числу ребер графа $F$. Это доказательство изложено в разделе~\ref{copies_lem_proof}.

Наконец, мы нашли явное представление $M'_{(R,H)}$ из
леммы~\ref{ext1}. Так как вывод этого представления опирается на
лемму~\ref{copies_lem}, то мы формулируем соответствующее
утверждение отдельно от леммы~\ref{ext1}.

\begin{lemma} \label{ext2}
В обозначениях леммы \ref{ext1} в качестве $M'_{(R,H)}$ можно
выбрать $M(k, l)$.
\end{lemma}

Лемма \ref{ext2} будет доказана в разделе~\ref{ext2_proof}.

\section{Доказательства лемм}
\label{lemmas_proofs}

Прежде чем перейти к доказательствам, введем вспомогательные
обозначения. Пусть $(R, H)$ --- произвольная сеть с $V(H)=\{z_1,
\ldots, z_d, y_1,\ldots, y_k\}$, $R = \{z_1, \ldots, z_d\}$ и
$e(R, H) = l$, а $\mathbf{v}^1, \ldots, \mathbf{v}^d\in V$ ---
произвольные вершины. В обозначениях из раздела~\ref{new_results}
положим
\begin{multline*}
    w_1^{\vec{x}}(1) = \left[n/2^d\right] + x_1, \ldots,
    w_{2^d-1}^{\vec{x}}(1) = \left[n/2^d\right] + x_{2^d-1},\\
    w_{2^d}^{\vec{x}}(1) = n - \left(2^d-1\right) \left[n/2^d\right] + x_{2^d}.
\end{multline*}
Рассмотрим произвольное $s\in\{1,\ldots,k\}$ и такие вершины
$\mathbf{w}^1,\ldots,\mathbf{w}^s\in V$, что $\{z_i, y_j\} \in
E(H) \Rightarrow \{\mathbf{v}^i, \mathbf{w}^j\} \in E(G)$ для
любых $i\in\{1,\ldots,d\}$, $j\in\{1,\ldots,s\}$ и $\{y_i, y_j\}
\in E(H) \Rightarrow \{\mathbf{w}^i, \mathbf{w}^j\} \in E(G)$ для
любых $i,j\in\{1,\ldots,s\}$. Для вершин
$\mathbf{v}^1,\ldots,\mathbf{v}^d,\mathbf{w}^1,\ldots,\mathbf{w}^{s-1}$
рассмотрим разбиение множества $\{1,\ldots,n\}$ на подмножества
$B_1,\ldots,B_{2^{d+s-1}}$ (см. раздел~\ref{new_results}),
мощности которых, как несложно видеть, зависят от чисел
$x_1,\ldots,x_{2^d}$, и положим
\[
   w_j=w_j^{\vec{x}}(s) = |B_j|,\quad
   B_j=\left\{r^j_1,\ldots,r^j_{w_j}\right\},\quad
   j\in\left\{1,\ldots,2^{d+s-1}\right\}.
\]
Вспомним, что вершина $\mathbf{w}^s$ есть вектор из $\{0, 1\}^n$. Посмотрим
на ее координаты. Для каждого $j\in\{1,\ldots,2^{d+s-1}\}$ обозначим
$u^{\vec{x}}_j(s)$ количество единиц среди чисел
$w^s_{r^j_1},\ldots,w^s_{r^j_{w_j}}$ (мы обозначили
$w^s_1,\ldots,w^s_n$ координаты вектора $\mathbf{w}^s$). Поскольку
$\mathbf{w}^s \in V$, ее скалярный квадрат равен $n/2$. Левую часть этого
равенства можно записать, очевидно, в виде суммы всех $u^{\vec{x}}_j(s)$,
$j\in\{1,\ldots,2^{d+s-1}\}$.
В силу определения
нашего графа $G$ вершина $\mathbf{w}^s$ соединена с какой-то вершиной
$\mathbf{v}$ из
$\mathbf{v}^1,\ldots,\mathbf{v}^d,\mathbf{w}^1,\ldots,\mathbf{w}^{s-1}$ тогда
и только тогда, когда их скалярное произведение равно $n/4$. Левая часть этого
равенства записывается в виде суммы $2^{d+s-2}$ величин $u^{\vec{x}}_j(s)$
(здесь индексы входящих в выражение величин суть индексы тех множеств $B_j$,
для которых координаты вершины $\mathbf{v}$ с номерами из $B_j$ равны $1$:
нулевые координаты в скалярном произведении не участвуют).
Предположим, что в графе $H$ среди вершин
$x_1,\ldots,x_d,y_1,\ldots,y_{s-1}$ только вершины $
   x_{l^1_1(s)},\ldots$, $x_{l^1_{a(s)}(s)},
   y_{l^2_1(s)},\ldots,y_{l^2_{b(s)}(s)}
$ соединены ребрами с вершиной $y_s$ (здесь $a(s)$ --- количество
вершин среди $x_1,\ldots,x_d$, соединенных с $y_s$, $b(s)$ ---
среди $y_1,\ldots,$ $y_{s-1}$). Обозначим $c(1, i)$, $i \in
\{l^1_1(s), \ldots, l^1_{a(s)}(s)\}$, и $c(2, i)$, $i \in
\{l^2_1(s), \ldots,$ $l^2_{b(s)}(s)\}$, последовательности индексов
переменных $u_j^{\vec{x}}(s)$, входящих в уравнения, соответствующие
наличию ребер между вершинами $y_s$ и $x_i$, $i \in \{l^1_1(s),
\ldots, l^1_{a(s)}(s)\}$, и между вершинами $y_s$ и $y_i$, $i \in
\{l^2_1(s), \ldots, l^2_{b(s)}(s)\}$, соответственно. Заметим, что
длины всех таких последовательностей совпадают и равны $m(s) = 2^{d+s-2}$.
Тогда для того, чтобы вершина
$\mathbf{w}^s$ была соединена с вершинами
$\mathbf{v}_{l^1_1(s)},\ldots,\mathbf{v}_{l^1_{a(s)}(s)}$,
$\mathbf{w}_{l^2_1(s)},\ldots,\mathbf{w}_{l^2_{b(s)}(s)}$, необходимо и
достаточно, чтобы были справедливы равенства и неравенства системы
\begin{equation}
\left\{
\begin{array}{c}
  u^{\vec{x}}_{c_1\left(1, l_1^1(s)\right)}(s)+\ldots+
  u^{\vec{x}}_{c_{m(s)}\left(1, l_1^1(s)\right)}(s)=n/4, \\
    \ldots \\
  u^{\vec{x}}_{c_1\left(1, l_{a(s)}^1(s)\right)}(s)+
  \ldots+u^{\vec{x}}_{c_{m(s)}\left(1, l_{a(s)}^1(s)\right)}(s)=n/4, \\

  u^{\vec{x}}_{c_1\left(2, l_1^2(s)\right)}(s)+\ldots+
  u^{\vec{x}}_{c_{m(s)}\left(2, l_1^2(s)\right)}(s)=n/4, \\
    \ldots \\
  u^{\vec{x}}_{c_1\left(2, l_{b(s)}^2(s)\right)}(s)+\ldots+u^{\vec{x}}_{c_{m(s)}\left(2, l_{b(s)}^2(s)\right)}(s)=
  n/4, \\
  u_{1}(s)+u_{2}(s)+\ldots+u_{2^{d+s-1}}(s)=n/2, \\
  \forall j\in \left\{1,2,3,\ldots,2^{d+s-1}\right\} \,\,\, 0\le u^{\vec{x}}_j(s)\le w^{\vec{x}}_j(s). \\
\end{array}
\right. \label{system}
\end{equation}
Очевидно,
\begin{equation}
    M_{(R, H)}^{\vec{x}} =
    \sum C_{w^{\vec{x}}_1(1)}^{u^{\vec{x}}_1(1)} \ldots
    C_{w^{\vec{x}}_{2^d}(1)}^{u^{\vec{x}}_{2^d}(1)} \cdot \ldots \cdot
    C_{w_1^{\vec{x}}(k)}^{u_1^{\vec{x}}(k)} \ldots
    C_{w_{2^{d+k-1}}^{\vec{x}}(k)}^{u_{2^{d+k-1}}^{\vec{x}}(k)},
\label{sum}
\end{equation}
где суммирование ведется по всем решениям
$(u^{\vec{x}}_1(1),\ldots,u^{\vec{x}}_{2^d}(1)),$ $\ldots,$
$(u^{\vec{x}}_1(k),$ $\ldots,$ $u^{\vec{x}}_{2^{d+k-1}}(k))$
систем~(\ref{system}) с $s=1$, $\ldots$, $s=k$ соответственно,
$w^{\vec{x}}_{2j-1}(s+1)=u^{\vec{x}}_j(s)$,
$w^{\vec{x}}_{2j}(s+1)=w^{\vec{x}}_j(s)-u^{\vec{x}}_j(s)$ при
любых $s\in\{1,\ldots,k-1\}$, $j\in\{1,\ldots,2^{d+s-1}\}$.

\subsection{Доказательство леммы \ref{ext1}}
\label{ext1_proof}
Без ограничения общности будем считать, что $f\gg n^{0.6}$.
Будем также предполагать, что $f(n) \ge n^{0.6}$ при всех $n$.

Мы начнем с идеи доказательства леммы, после чего
приведем полное строгое доказательство.

\paragraph{Идея доказательства}

В разделе \ref{lemmas_proofs} мы свели задачу нахождения числа расширений
$M_{(R, H)}^{\vec{x}}$ к поиску асимптотики суммы~(\ref{sum})
по всем решениям $k$ систем уравнений вида~(\ref{system}). Требуется доказать, что эта асимптотика
не зависит от $x_1, \ldots, x_{2^d}$ и равномерна по ним при $|x_i| \le
f(n)$, $i \in \{1, \ldots, 2^d\}$.

Нетрудно видеть, что максимум суммы~(\ref{sum}) достигается, когда
$u_j^{\vec{x}}(s) = $ $[w_j^{\vec{x}}(s)/2] + O(1)$, $n \to \infty$. Более того,
далее будет показано, что если в сумме~(\ref{sum}) проводить суммирование не по
всем $u_j^{\vec{x}}(s)$, удовлетворяющим соответствующим системам
уравнений~(\ref{system}), а лишь по таким, что $|u_j^{\vec{x}}(s) -
[w_j^{\vec{x}}(s)/2]| \le n^{0.6}$, то асимптотика суммы не поменяется,
причем сохранится и равномерность по $x_1, \ldots, x_{2^d}$ (при вышеупомянутых
ограничениях). Для такой ``укороченной'' суммы уже гораздо легче, вводя для
удобства новые переменные и обозначения, доказать равномерность ее асимптотики.

\paragraph{Доказательство}

В силу определения чисел $x_1,\ldots,x_{2^d}$ существуют такая
константа $c>0$ и такие целые числа $a_1,\ldots,a_d$, что
$|a_j|\le c$ для всех $j\in\{1,\ldots,d\}$ и вектор
$(x_1,\ldots,x_{2^d})$ является решением системы~(\ref{system}), в
которой $s=1$, $a(1)=d$ и правая часть заменена на столбец
$(a_1,\ldots,a_d,0)^T$. Докажем, что существует такая константа
$C>0$, не зависящая от $(x_1,\ldots,x_{2^d})$, что для всех $j\in \{1,\ldots,2^d\}$ найдутся
числа $y_j=y_j(n)\in\mathbb{Z}$ и $r_j\in\mathbb{Z}$,
\begin{equation}
   |r_j|\le C,
\label{const_c}
\end{equation}
для которых
\begin{equation}
   x_j = 2^k y_j + r_j,
\label{from_x_to_y}
\end{equation}
а вектор $(y_1,\ldots,y_{2^d})$ является решением
системы~(\ref{system}), в которой $s=1$, $a(1)=d$ и правая часть
заменена на столбец $(0,\ldots,0)^T$. Обозначим последнюю систему
следующим образом: $A\mathbf{y}=\mathbf{0}$. Подставив вместо
$\mathbf{y}$ вектор $([x_1/2^k],\ldots,[x_{2^d}/2^k])$, получим в
правой части некоторый вектор $(b_1,\ldots,$ $b_{d+1})^T$, абсолютное
значение каждого элемента которого не превосходит некоторой
константы $\tilde c>0$. Заметим, что для любого
$i\in\{1,\ldots,d+1\}$ в $i$-ой строке матрицы
$A=(a_{i,j})_{d+1}^{2^d}$ содержится ненулевой коэффициент: $a_{i,2^{d}-2^{d-i}}=1$ (если $i\in\{1,\ldots,d\}$) и
$a_{d+1,2^d}=1$, при этом соответствующие коэффициенты в остальных
строках равны нулю: $a_{j,2^{d}-2^{d-i}}=0$ (если
$i\in\{1,\ldots,d\}$, $j\in\{1,\ldots,i-1,i+1,\ldots,d\}$) и
$a_{j,2^d}=0$ (если $j\in\{1,\ldots,d\}$). Следовательно, положив
\[
   y_{2^{d}-2^{d-i}}=\left[\frac{x_{2^d-2^{d-i}}}{2^{k}}\right]-b_i \text{ при } i\in\{1,\ldots,d\},
\]
\[
   y_{2^d}=\left[\frac{x_{2^d}}{2^k}\right] + \sum_{i=1}^d b_i - b_{d+1},
\]
\[
   y_j=\left[\frac{x_j}{2^k}\right] \text{ при }
   j\in\left\{1,\ldots,2^d\right\}\setminus\left\{2^d-2^{d-1},2^d-2^{d-2},
   \ldots,2^d-1,2^d\right\},
\]
мы получим искомый вектор $\mathbf{y}$.

Введем для удобства новые обозначения:
\[
   w_j'(1)=\left[\frac{n}{2^d}\right] \text{ при }
   j\in\left\{1,\ldots,2^d-1\right\},\quad
   w'_{2^d}(1)=n-\left(2^d-1\right)\left[\frac{n}{2^d}\right],
\]
\begin{equation}
   t_j(s)=u^{\vec{x}}_j(s)-\left[\frac{w^{\vec{x}}_j(s)-\varepsilon_j(s)}{2}\right],
\label{t}
\end{equation}
где $\varepsilon_j(1)=r_j$ при $j\in\{1,\ldots,2^d\}$, а при
$s\in\{2,\ldots,k\}$
\[
    \varepsilon_j(s) =
    \begin{cases}
        r_q,  & \text{если $j = 2^{s-1} q$, $q \in \{1, \ldots, 2^d\}$;} \\
        0,  & \text{иначе.}
    \end{cases}
\]
Далее $u_j'(s)$ при $s\in\{1,\ldots,k\}$,
$j\in\{1,\ldots,2^{d+s-1}\}$ и $w_j'(s)$ при $s\in\{2,\ldots,k\}$,
$j\in\{1,\ldots,2^{d+s-1}\}$ определяются из следующих
рекуррентных соотношений
\begin{equation}
   u_j'(s) = \left[w_j'(s)/2\right] + t_j(s),
\label{recurcy_first}
\end{equation}
\begin{equation}
    w_{2 j - 1}'(s) = u_j'(s-1), \quad
    w_{2 j}'(s) = w_j'(s-1) - u_j'(s-1).
\label{recurcy_second}
\end{equation}
Очевидно, $w_j^{\vec{x}}(1)=w_j'(1)+2^k y_j+r_j$ и
$u_j^{\vec{x}}(1)=[w_j'(1)/2]+2^{k-1}y_j+t_j(1)$ при
$j\in\{1,\ldots,2^d\}$. Получим по индукции аналогичные выражения
для $w_j^{\vec{x}}(s)$ и $u_j^{\vec{x}}(s)$ при
$s\in\{2,\ldots,k\}$, $j\in\{1,\ldots,2^{d+s-1}\}$. Пусть
$s\in\{2,\ldots,k\}$. Предположим, что для любого
$j\in\{1,\ldots,2^{d+(s-1)-1}\}$ справедливы равенства
\begin{equation}
\begin{gathered}
   w^{\vec{x}}_j(s-1)=w_j'(s-1) + 2^{k-(s-1)+1} \delta_j(s-1) +
   \varepsilon_j(s-1),\\
   u^{\vec{x}}_j(s-1)=[w_j'(s-1)/2] +
   2^{k-(s-1)} \delta_j(s-1) + t_j(s-1),
\label{previous_step}
\end{gathered}
\end{equation}
где
\[
    \delta_j(s) =
    y_q, \text{ если $j \in \left\{2^{s-1} (q-1) + 1, \ldots, 2^{s-1} q\right\}$, $q \in \left\{1, \ldots, 2^d\right\}$.}
\]
Тогда в
силу~(\ref{recurcy_first}),~(\ref{recurcy_second})~и~(\ref{previous_step})
для любого $j\in\{1,\ldots,2^{d+(s-1)-1}\}$ имеем
\[
   w_{2j-1}^{\vec{x}}(s)=u_j^{\vec{x}}(s-1)=u'_j(s-1)+2^{k-s+1}\delta_j(s-1)=
   w'_{2j-1}(s)+2^{k-s+1}\delta_{2j-1}(s),
\]
\begin{multline*}
   w_{2j}^{\vec{x}}(s)=w_j^{\vec{x}}(s-1)-u_j^{\vec{x}}(s-1)=
   w'_j(s-1)-u'_j(s-1)+2^{k-s+1}\delta_j(s-1)
   +\\
   \varepsilon_j(s-1)=
   w'_{2j}(s)+2^{k-s+1}\delta_{2j}(s)+\varepsilon_{2j}(s).
\end{multline*}
Окончательно, для любого $j\in\{1,\ldots,2^{d+s-1}\}$ выполнено
\begin{equation}
   w_{j}^{\vec{x}}(s)=w'_{j}(s)+2^{k-s+1}\delta_{j}(s)+\varepsilon_{j}(s).
\label{from_w_to_w}
\end{equation}
Следовательно, в силу~(\ref{t})
\begin{equation}
   u^{\vec{x}}_j(s)=t_j(s)+\left[\frac{w^{\vec{x}}_j(s)-\varepsilon_j(s)}{2}\right]=
   \left[\frac{w_j'}{2}\right]+2^{k-s}\delta_j(s)+t_j(s).
\label{u}
\end{equation}
Таким образом,
\[
    M_{(R, H)}^{\vec{x}} =
    \sum C_{w_1'(1) + 2^k \delta_1(1) + \varepsilon_1(1)}^
    {\left[w_1'(1)/2\right] + 2^{k-1} \delta_1(1) + t_1(1)} \ldots
    C_{w_{2^{d}}'(1) + 2^k \delta_{2^d}(1) + \varepsilon_{2^d}(1)}^
    {\left[w_{2^{d}}'(1)/2\right]+2^{k-1}\delta_{2^d}(1)+t_{2^d}(1)}
    \ldots\times
\]
\begin{equation}
    C_{w_1'(k) + 2 \delta_1(k) + \varepsilon_1(k)}^
    {\left[w_1'(k)/2\right] + \delta_1(k) + t_1(k)} \ldots
    C_{w_{2^{d+k-1}}'(k) + 2 \delta_{2^{d+k-1}}(k) + \varepsilon_{2^{d+k-1}}(k)}^
    {\left[w_{2^{d+k-1}}'(k)/2\right] + \delta_{2^{d+k-1}}(k) + t_{2^{d+k-1}}(k)},
\label{from_u_to_t_1}
\end{equation}
где
суммирование ведется по всем таким наборам
$(t_1(1),\ldots,t_{2^d}(1),\ldots,t_1(k),$ $\ldots,$ $t_{2^{d+k-1}}(k))$,
что числа $u_j^{\vec{x}}(s)$, $s\in\{1,\ldots,k\}$,
$j\in\{1,\ldots,2^{d+s-1}\}$, определяемые равенством~(\ref{t}),
являются решением систем~(\ref{system}) при $s\in\{1,\ldots,k\}$.

Докажем, что суммирование в~(\ref{from_u_to_t_1}) можно осуществлять
по некоторому такому множеству $T=T(n)$ наборов
$(t_1(1),\ldots,t_{2^d}(1),\ldots,t_1(k),\ldots,$ $t_{2^{d+k-1}}(k))$,
не зависящему от $x_1,\ldots,x_{2^d}$, что равенство
в~(\ref{from_u_to_t_1}) заменится на асимптотическое равенство, и
для всех
$(t_1(1),\ldots,t_{2^d}(1),\ldots,t_1(k),\ldots,$ $t_{2^{d+k-1}}(k))\in
T$ будет выполнено $|t_j(s)|\le n^{0.6}$ при
$s\in\{1,\ldots,k\}$, $j\in\{1,\ldots,2^{d+s-1}\}$.

Заметим, что в силу определения вектора $(y_1,\ldots,y_{2^d})$
вектор $(\delta_1(1),\ldots,$ $\delta_{2^d}(1))=(y_1,\ldots,y_{2^d})$
является решением системы~(\ref{system}), в которой $s=1$,
$a(1)=d$ и правая часть заменена на столбец $(0,\ldots,0)^T$. Произведем замену переменных $u^{\vec{x}}_j(1)$, $j \in \{1, \ldots, 2^d\}$, в системе~(\ref{system}) в соответствии с (\ref{u}). Тогда коэффициенты системы при $s=1$ без ограничений, записанных в ее последней строке, не зависят от $\vec{x}$. Следовательно, и ее решения также не зависят от $\vec{x}$. Более того, при достаточно
больших $n$ любое решение $(t_1(1),\ldots,t_{2^d}(1))$
системы~(\ref{system}) при $s=1$ без ограничений, записанных в ее
последней строке, удовлетворяющее неравенствам $|t_j(1)|\le
n^{0.6}$, $j\in\{1,\ldots,2^d\}$, удовлетворяет этим ограничениям,
так как $w_j^{\vec{x}}(1)\ge[n/2^d]-f(n)$ при всех
$j\in\{1,\ldots,2^d\}$. Следовательно, множество всех наборов
$(t_1(1),\ldots,t_{2^d}(1))$, удовлетворяющих неравенствам
$|t_j(1)|\le n^{0.6}$, $j\in\{1,\ldots,2^d\}$, и определяемых
равенством~(\ref{t}), совпадает с множеством всех наборов
$(t_1(1),\ldots,t_{2^d}(1))$, удовлетворяющих неравенствам
$|t_j(1)|\le n^{0.6}$, $j\in\{1,\ldots,2^d\}$, и являющихся
решениями системы~(\ref{system}) без ограничений, записанных в ее
последней строке. Обозначим $\tilde T(1)$ множество всех таких
наборов. Заметим, наконец, что в
силу~(\ref{from_x_to_y}),~(\ref{recurcy_first}),~(\ref{recurcy_second})~и~(\ref{from_w_to_w})
для любого набора $(t_1(1),\ldots,t_{2^d}(1))\in \tilde T(1)$ выполнены
неравенства
\[
   w_j^{\vec{x}}(2) \ge
   w_j'(2) - C_1 f(n) \ge \left[w_j'(1)/2\right] - n^{0.6} - C_1 f(n) \ge \left[n/2^{d+1}\right]-\tilde C(2) f(n)
\]
при всех $j\in\left\{1,\ldots,2^{d+1}\right\}$, где $C_1$ и $\tilde C(2)$ --- некоторые положительные константы.

Пусть $s\in\{2,\ldots,k\}$ и множество $\tilde T(s-1)$ наборов
$(t_1(1),\ldots,t_{2^d}(1),\ldots,$ $t_1(s-1),\ldots,t_{2^{d+s-2}}(s-1))$
задано. Рассмотрим произвольный набор
\[
   \mathbf{t}=(t_1(1),\ldots,t_{2^d}(1),\ldots,t_1(s-1),\ldots,t_{2^{d+s-2}}(s-1))\in\tilde
   T(s-1).
\]
Пусть для всех $j\in\{1,\ldots,2^{d+s-1}\}$ выполнено
\begin{equation}
   w_j^{\vec{x}}(s)\ge \left[\frac{n}{2^{d+s-1}}\right]-
   \tilde C(s) f(n),
\label{w_bound}
\end{equation}
где $\tilde C(s) = \const > 0$.
По аналогии со случаем $s=1$ в силу определения вектора
$(y_1,\ldots,y_{2^d})$ вектор $(\delta_1(s),\ldots,\delta_{2^{d+s-1}}(s))$ является решением системы~(\ref{system}), в которой $a(s)=d$, $b(s)=s-1$ и правая часть заменена на столбец $(0,\ldots,0)^T$. Снова сделаем замену переменных $u^{\vec{x}}_j(s)$, $j \in \{1, \ldots, 2^{d+s-1}\}$, в соответствии с (\ref{u}). Поскольку коэффициенты системы без ограничений, записанных в ее последней строке, не зависят от $\vec{x}$ и являются функциями от переменных $t_j(i)$, $j \in \{1, \ldots, 2^{d+i-1}\}$, $i \in \{1, \ldots, s-1\}$, которые также не зависят от $\vec{x}$, то и решения такой системы не зависят от $\vec{x}$, и при достаточно больших $n$ любое решение
$(t_1(s),\ldots,t_{2^{d+s-1}}(s))$ системы~(\ref{system}) без
ограничений, записанных в ее последней строке, удовлетворяющее
неравенствам $|t_j(s)|\le n^{0.6}$, $j\in\{1,\ldots,2^{d+s-1}\}$,
удовлетворяет этим ограничениям в силу предположения~(\ref{w_bound}). Следовательно, множество всех наборов $(t_1(s),\ldots,t_{2^{d+s-1}}(s))$, удовлетворяющих неравенствам $|t_j(s)|\le n^{0.6}$, $j\in\{1,\ldots,2^{d+s-1}\}$, и определяемых равенством~(\ref{t}), совпадает с множеством всех наборов $(t_1(s),\ldots,t_{2^{d+s-1}}(1))$, удовлетворяющих
неравенствам $|t_j(s)|\le n^{0.6}$, $j\in\{1,\ldots,2^{d+s-1}\}$, и
являющихся решениями системы~(\ref{system}) без ограничений,
записанных в ее последней строке. Обозначим $\tilde T(\mathbf{t})$
множество всех таких наборов. При $s<k$,
$j\in\{1,\ldots,2^{d+s}\}$ и любом
$(t_1(s),\ldots,t_{2^{d+s-1}}(s))\in\tilde T(\mathbf{t})$ в
силу~(\ref{from_x_to_y}),~(\ref{recurcy_first}),~(\ref{recurcy_second}),~(\ref{from_w_to_w})
и предположения~(\ref{w_bound}) имеем
\begin{multline*}
   w_j^{\vec{x}}(s+1)\ge
   w_j'(s+1) - C_1 f(n) \ge
   \left[\frac{w_j'(s)}{2}\right]-n^{0.6}-C_1 f(n)\ge \\
   \left[\frac{w_j^{\vec{x}}(s)-2^{k-s+1}\delta_j(s)-
   \varepsilon_j(s)}{2}\right]-C_2 f(n) \ge
   \left[\frac{n}{2^{d+s}}\right]-\tilde C(s+1) f(n),
\end{multline*}
где $C_1, C_2, \tilde C(s+1)$ --- положительные константы.
Положим
\[
 \tilde T(s)=\bigcup_{\mathbf{t}=(t_1(1),\ldots,t_{2^d}(1),\ldots,t_1(s-1),\ldots,t_{2^{d+s-2}}(s-1))\in
 \tilde T(s-1)} \{\mathbf{t}\}\times\tilde T(\mathbf{t}).
\]
Множество $T$ определим следующим образом: $T=\tilde T(k)$.
Осталось доказать, что равенство в~(\ref{from_u_to_t_1}) заменится
на асимптотическое равенство при суммировании по всем наборам из
$T$.

По аналогии с доказательством существования чисел
$y_1,\ldots,y_{2^d}$, для которых
выполнено~(\ref{const_c})~и~(\ref{from_x_to_y}), можно доказать,
что для любых $x_1, \ldots, x_{2^d}$ существует такое число
$c=\const>0$ и набор $(t_1(1),\ldots, t_{2^{d}}(1), \ldots,
t_1(k), \ldots, $ $t_{2^{d+k-1}}(k))\in T$, что выполнены неравенства $t_j(s)\le c$ для
всех $s\in\{1,\ldots,k\}$, $j\in\{1,\ldots,$ $2^{d+s-1}\}$.
Соответствующий этому набору элемент суммы из правой части
равенства~(\ref{from_u_to_t_1}), очевидно, оценивается снизу
величиной $2^{k n}/n^{\tilde c}$ для некоторой константы $\tilde
c>0$.

Заметим, что по формуле Стирлинга для любых натуральных чисел $z$
и $d$ выполнено неравенство $C_z^{[z/2] + d}/C_z^{[z/2]} \le
e^{-2d^2/z} (1 + \varepsilon(d,z))$, где
$\varepsilon_d(z)\rightarrow 0$ при $z\rightarrow\infty$ и
$d/z\rightarrow 0$. При фиксированном $z$ величина
$C_z^{[z/2] + d}/C_z^{[z/2]}$ убывает при увеличении $d$, так что
$C_z^{[z/2] + d}/C_z^{[z/2]} \le e^{-2 z^{0.2}} (1 +
\varepsilon(z))$ при $|d| > z^{0.6}$, где
$\varepsilon(z)\rightarrow 0$ при $z\rightarrow\infty$. Поэтому
слагаемые в~(\ref{from_u_to_t_1}), соответствующие наборам
$((t_1(1),\ldots, t_{2^{d}}(1), \ldots, t_1(k), \ldots,
t_{2^{d+k-1}}(k)))$, в которых хотя бы одно из $t_j(i)$
удовлетворяет неравенству $\left|t_j(i)\right| > n^{0.6}$, дают в
совокупности
\[
 n^{k2^{d+k-1}}e^{- \Omega\left(n^{0.2}\right)}2^{kn}=
 o\left(2^{k n}/n^{\tilde c}\right).
\]
Здесь $g(n) = \Omega(h(n))$ означает, что при достаточно больших $n$ для некоторого $C>0$ выполнено $g(n) \ge C h(n)$.
Окончательно получаем
\begin{multline}
    M_{(R, H)}^{\vec{x}} =
    \sum\limits_{(t_1(1),\ldots,t_{2^d}(1),\ldots,t_1(k),\ldots,t_{2^{d+k-1}}(k))\in T}
    C_{w_1'(1) + 2^k \delta_1(1) + \varepsilon_1(1)}^
    {\left[w_1'(1)/2\right] + 2^{k-1} \delta_1(1) + t_1(1)} \ldots \times \\
    C_{w_{2^{d}}'(1) + 2^k \delta_{2^d}(1) + \varepsilon_{2^d}(1)}^
    {\left[w_{2^{d}}'(1)/2\right]+2^{k-1}\delta_{2^d}(1)+t_{2^d}(1)}
    \ldots
    C_{w_1'(k) + 2 \delta_1(k) + \varepsilon_1(k)}^
    {\left[w_1'(k)/2\right] + \delta_1(k) + t_1(k)} \ldots \times \\
    C_{w_{2^{d+k-1}}'(k) + 2 \delta_{2^{d+k-1}}(k) + \varepsilon_{2^{d+k-1}}(k)}^
    {\left[w_{2^{d+k-1}}'(i)/2\right] + \delta_{2^{d+k-1}}(k) +
    t_{2^{d+k-1}}(k)}\left(1+o\left(\frac{1}{n}\right)\right)
\label{from_u_to_t}
\end{multline}
равномерно по всем $x_1,\ldots,x_{2^d} \in[-f(n), f(n)]$.

Осталось теперь доказать, что асимптотика каждого члена суммы по множеству $T$
равномерна по $x_1,\ldots,x_{2^d} \in[-f(n), f(n)]$,
$(t_1(1), \ldots, t_{2^d}(1), \ldots,$ $t_1(k), \ldots,$ $t_{2^{d+k-1}}(k))\in T$ и $s \in \left\{1, \ldots, k\right\}$. Покажем, что
\begin{multline*}
    C_{w_1'(s) + 2^{k-s+1} \delta_1(s) + \varepsilon_1(s)}^
    {\left[w_1'(s)/2\right] + 2^{k-s} \delta_1(s) + t_1(s)} \ldots
    C_{w_{2^{d+s-1}}'(s) + 2^{k-s+1} \delta_{2^{d+s-1}}(s) + \varepsilon_{2^{d+s-1}}(s)}^
    {\left[w_{2^{d+s-1}}'(i)/2\right] + 2^{k-s} \delta_{2^{d+s-1}}(s) +
    t_{2^{d+s-1}}(s)}\sim \\
    C_{w_1'(s)}^{\left[w_1'(s)/2\right] + t_1(s)} \ldots C_{w_{2^{d+s-1}}'(s)}^
    {\left[w_{2^{d+s-1}}'(s)/2\right]+ t_{2^{d+s-1}}(s)}
\end{multline*}
равномерно по всем $x_1,\ldots,x_{2^d} \in[-f(n), f(n)]$,
$(t_1(1), \ldots, t_{2^d}(1), \ldots, t_1(k),$ $\ldots,$ $t_{2^{d+k-1}}(k))\in T$ и $s \in \left\{1, \ldots, k\right\}$.

Избавимся, в первую очередь, от $\varepsilon_j(s)$ в биномиальных
коэффициентах. Заметим, что в силу~(\ref{w_bound}) найдется такое $c>0$,
что для всех $x_1,\ldots,x_{2^d} \in[-f(n), f(n)]$, $(t_1(1), \ldots,
t_{2^d}(1), \ldots, t_1(k), \ldots, t_{2^{d+k-1}}(k))\in T$, $s
\in \left\{1, \ldots, k\right\}$ и $j\in\{1,\ldots,2^{d+s-1}\}$
выполнено $cn \le w_j^{\vec{x}}(s) \le n$. В этой связи по формуле
Стирлинга найдется такая константа $\tilde c>0$, что для всех тех
же наборов переменных
\[
    \left|\left(C_{w_j'(s) + 2^{k-s+1} \delta_j(s) + \varepsilon_j(s)}^
    {\left[w_j'(s)/2\right] + 2^{k-s} \delta_j(s) + t_j(s)}\right)\bigg/
    \left(C_{w_j'(s) + 2^{k-s+1} \delta_j(s)}^{\left[w_j'(s)/2\right] +
    2^{k-s} \delta_j(s) + t_j(s)}\right) -
    2^{\varepsilon_j(s)}\right|\le\tilde c\frac{f(n)}{n}.
\]
Более того, по формуле Стирлинга существует такое $M>0$, что для
всех $w\in[cn,n]\cap\mathbb{N}$ имеем
\[
    \left|\frac{w!}{\sqrt{2 \pi w} \frac{w^w}{e^w}}-1\right|\le
    \frac{M}{w}.
\]
Разложение Тейлора дает существование такого $c_1>0$, что для всех
$s\in\{1,\ldots,k\}$, $w\in[cn,n]$, $|y|\le f(n)$
\[
    \left|\frac{\sqrt{w}}{\sqrt{w + 2^{k - s + 1} y}}
    -1\right|\le c_1\left|\frac{y}{w}\right|\le \frac{c_1}{c}\cdot\frac{f(n)}{n}.
\]
Наконец, в силу определения величин $x_1,\ldots,x_{2^d}$ и
$y_1,\ldots,y_{2^d}$ справедливы равенства
$\sum_{j=1}^{2^d}x_j=\sum_{j=1}^{2^d}y_j=0$, а, следовательно,
$\sum_{j=1}^{2^d}r_j=0$ в силу~$(\ref{from_x_to_y})$ и
$\sum_{j=1}^{2^{d+s-1}}\delta_j(s)=\sum_{j=1}^{2^{d+s-1}}\varepsilon_j(s)=0$ для любого
$s\in\{1,\ldots,k\}$ в силу определения величин $\delta_j(s)$ и $\varepsilon_j(s)$.

Поэтому
\begin{multline}\label{first}
    C_{w_1'(s) + 2^{k-s+1} \delta_1(s) + \varepsilon_1(s)}^
    {\left[w_1'(s)/2\right] + 2^{k-s} \delta_1(s) + t_1(s)} \ldots
    C_{w_{2^{d+s-1}}'(s) + 2^{k-s+1} \delta_{2^{d+s-1}}(s) + \varepsilon_{2^{d+s-1}}(s)}^
    {\left[w_{2^{d+s-1}}'(s)/2\right] + 2^{k-s} \delta_{2^{d+s-1}}(s) + t_{2^{d+s-1}}(s)} = \\
    2^{\sum_{j=1}^{2^{d+s-1}} \varepsilon_j(s)} \times \\
    C_{w_1'(s) + 2^{k-s+1} \delta_1(s)}^
    {\left[w_1'(s)/2\right] + 2^{k-s} \delta_1(s) + t_1(s)} \ldots
    C_{w_{2^{d+s-1}}'(s) + 2^{k-s+1} \delta_{2^{d+s-1}}(s)}^
    {\left[w_{2^{d+s-1}}'(s)/2\right] + 2^{k-s} \delta_{2^{d+s-1}}(s) + t_{2^{d+s-1}}(s)} \times \\
    \left(1 + O\left(\frac{f(n)}{n}\right)\right) = \\
    \left(\sqrt{\frac{2}{\pi}}\right)^{2^{d+s-1}}
    \frac{2^{w_1'(s) + 2^{k - s + 1} \delta_1(s)}\ldots
    2^{w_{2^{d+s-1}}'(s) + 2^{k - s + 1} \delta_{2^{d+s-1}}(s)}}
    {\sqrt{\left(w_1'(s) + 2^{k - s + 1} \delta_1(s)\right)\ldots
    \left(w_{2^{d+s-1}}'(s) + 2^{k - s + 1} \delta_{2^{d+s-1}}(s)\right)}}\times\\
    \frac{C_{w_1'(s) + 2^{k-s+1} \delta_1(s)}^{\left[w_1'(s)/2\right] + 2^{k-s} \delta_1(s) + t_1(s)}}
    {C_{w_1'(s) + 2^{k-s+1} \delta_1(s)}^{\left[w_1'(s)/2\right] + 2^{k-s} \delta_1(s)}}\ldots
    \frac{C_{w_{2^{d+s-1}}'(s) + 2^{k-s+1} \delta_{2^{d+s-1}}(s)}^
    {\left[w_{2^{d+s-1}}'(s)/2\right] + 2^{k-s} \delta_{2^{d+s-1}}(s) + t_{2^{d+s-1}}(s)}}
    {C_{w_{2^{d+s-1}}'(s) + 2^{k-s+1} \delta_{2^{d+s-1}}(s)}^
    {\left[w_{2^{d+s-1}}'(s)/2\right] + 2^{k-s} \delta_{2^{d+s-1}}(s)}} \times \\
    \left(1+O\left(\frac{f(n)}{n}\right)\right) =
    \left(\sqrt{\frac{2}{\pi}}\right)^{2^{d+s-1}}
    \frac{2^{w_1'(s)}\ldots 2^{w_{2^{d+s-1}}'(s)}}{\sqrt{w_1'(s) \ldots w_{2^{d+s-1}}'(s)}}\times\\
    \frac{C_{w_1'(s) + 2^{k-s+1} \delta_1(s)}^{\left[w_1'(s)/2\right] + 2^{k-s} \delta_1(s) + t_1(s)}}
    {C_{w_1'(s) + 2^{k-s+1} \delta_1(s)}^{\left[w_1'(s)/2\right] + 2^{k-s} \delta_1(s)}}\ldots
    \frac{C_{w_{2^{d+s-1}}'(s) + 2^{k-s+1} \delta_{2^{d+s-1}}(s)}^
    {\left[w_{2^{d+s-1}}'(s)/2\right] + 2^{k-s} \delta_{2^{d+s-1}}(s) + t_{2^{d+s-1}}(s)}}
    {C_{w_{2^{d+s-1}}'(s) + 2^{k-s+1} \delta_{2^{d+s-1}}(s)}^
    {\left[w_{2^{d+s-1}}'(s)/2\right] + 2^{k-s} \delta_{2^{d+s-1}}(s)}} \times \\
    \left(1+O\left(\frac{f(n)}{n}\right)\right)
\end{multline}
равномерно по всем $x_1,\ldots,x_{2^d} \in[-f(n), f(n)]$,
$(t_1(1), \ldots, t_{2^d}(1), \ldots, t_1(k),$ $\ldots,$ $t_{2^{d+k-1}}(k))\in T$ и $s \in \left\{1, \ldots, k\right\}$.

Для любых $x_1,\ldots,x_{2^d} \in[-f(n), f(n)]$, $s \in \left\{1,
\ldots, k\right\}$ и $j\in\{1,\ldots,2^{d+s-1}\}$ обозначим
$W^{\vec{x}}_j(s)$ множество всех троек
$(w'_j(s),2^{k-s}\delta_j(s),t_j(s))$, которые возникают, когда
набор $(t_1(1),\ldots,t_{2^d}(1),\ldots,$
$t_1(k),\ldots,t_{2^{d+k-1}}(k))$ пробегает все множество $T$.
Равномерно по всем $x_1,\ldots,x_{2^d} \in[-f(n), f(n)]$,
$s \in \left\{1, \ldots, k\right\}$,
$j\in\{1,\ldots,2^{d+s-1}\}$ и $(w,y,t)\in W^{\vec{x}}_j(s)$
\[
    \frac{C_{w + 2 y}^{[w/2] + y + t}}{C_{w + 2 y}^{[w/2] + y}} =
    e^{-2 t^2/w} \left(1 + O\left(f(n) n^{-0.8}\right)\right).
\]
Таким образом, в силу (\ref{first}) равномерно по всем
$x_1,\ldots,x_{2^d} \in[-f(n), f(n)]$, $(t_1(1), \ldots,
t_{2^d}(1), \ldots, t_1(k),$ $\ldots, t_{2^{d+k-1}}(k))\in T$ и $s
\in \left\{1, \ldots, k\right\}$
\begin{multline*}
    C_{w_1'(s) + 2^{k-s+1} \delta_1(s) + \varepsilon_1(s)}^
    {\left[w_1'(s)/2\right] + 2^{k-s} \delta_1(s) + t_1(s)} \ldots
    C_{w_{2^{d+s-1}}'(s) + 2^{k-s+1} \delta_{2^{d+s-1}}(s) + \varepsilon_{2^{d+s-1}}(s)}^
    {\left[w_{2^{d+s-1}}'(s)/2\right] + 2^{k-s} \delta_{2^{d+s-1}}(s) + t_{2^{d+s-1}}(s)} = \\
    \left(\sqrt{\frac{2}{\pi}}\right)^{2^{d+s-1}} \frac{2^{w_1'(s)}\ldots 2^{w_{2^{d+s-1}}'(s)}}
    {\sqrt{w_1'(s) \ldots w_{2^{d+s-1}}'(s)}} \times \\
    \frac{C_{w_1'(s)}^{\left[w_1'(s)/2\right] + t_1(s)}}{C_{w_1'(s)}^{\left[w_1'(s)/2\right]}}
    \ldots \frac{C_{w_{2^{d+s-1}}'(s)}^{\left[w_{2^{d+s-1}}'(s)/2\right] + t_{2^{d+s-1}}(s)}}
    {C_{w_{2^{d+s-1}}'(s)}^{\left[w_{2^{d+s-1}}'(s)/2\right]}}\left(1+O\left(f(n) n^{-0.8}\right)\right)= \\
    C_{w_1'(s)}^{\left[w_1'(s)/2\right] + t_1(s)}
    \ldots C_{w_{2^{d+s-1}}'(s)}^{\left[w_{2^{d+s-1}}'(s)/2\right] + t_{2^{d+s-1}}(s)}\left(1+O\left(f(n) n^{-0.8}\right)\right),
\end{multline*}
что, в свою очередь, доказывает лемму, так как множество $T$ и все
величины в правой части последнего равенства не зависят от
$\vec{x}$.

\subsection{Доказательство леммы \ref{copies_lem}}
\label{copies_lem_proof}
Как можно судить по доказательству предыдущей леммы, утверждение
этой леммы напрямую доказывать было бы тяжело, если вообще возможно. К счастью,
проблема легко решается, если доказывать утверждение сразу для всех связных
фиксированных графов.

Обозначим $\mathcal{M}_{m,l}$ множество всех связных графов с $m$
вершинами и $l$ ребрами.

Будем доказывать утверждение для каждого $m$ индукцией по числу
ребер. При $l = m - 1$ первую вершину (корень дерева)
выбираем $N$
способами, а каждую последующую (соединенную ребром с ровно одной
из ранее выбранных) $N_1 - O(1)$ способами (напомним, что $N_1$ ---
степень вершины рассматриваемого дистанционного графа). Поэтому
для любого $F\in\mathcal{M}_{m,m-1}$ имеем
\[
    M_F \sim N N_1^{m - 1}=N^m \left(\frac{N_1}{N}\right)^{m-1}.
\]

Пусть лемма верна для $l \le L - 1$, докажем ее для $l = L$, где
$L\in\{m,m+1,\ldots\}$. Пусть, кроме того,
$F\in\mathcal{M}_{m,L-1}$ и вершины $z_1$ и $z_2$ графа $F$ не
соединены ребром.

Определим сеть $(R, H)$ следующим образом: $H=F$, $R=\{z_1,z_2\}$. Рассмотрим произвольные вершины $\mathbf{v}^1=
(v_1^1,\ldots,v_n^1),\mathbf{v}^2=(v_1^2,\ldots,v_n^2) \in V$. Тогда в наших обозначениях $\vec{x}=(x_1,x_2,x_3,x_4)=(x,-x,-x,x)$ для некоторого $x\in\mathbb{Z}$. Поэтому
\begin{equation}
    \Sigma(x):=M_{(R, H)}^{\vec{x}} =
    \sum C_{w^{\vec{x}}_1(1)}^{u^{\vec{x}}_1(1)} \ldots
    C_{w^{\vec{x}}_{4}(1)}^{u^{\vec{x}}_{4}(1)} \cdot \ldots \cdot
    C_{w_1^{\vec{x}}(m-2)}^{u_1^{\vec{x}}(m-2)} \ldots
    C_{w_{2^{m-1}}^{\vec{x}}(m-2)}^{u_{2^{m-1}}^{\vec{x}}(m-2)},
\label{sum2}
\end{equation}
где суммирование ведется по всем решениям
$(u^{\vec{x}}_1(1),\ldots,u^{\vec{x}}_{4}(1)),$ $\ldots,$
$(u^{\vec{x}}_1(m-2),\ldots,u^{\vec{x}}_{2^{m-1}}(m-2))$
систем~(\ref{system}) с $s=1$, $\ldots$, $s=m-2$ соответственно.
Кроме того,
\[
   M_F=N\sum_{x=-n/4}^{n/4-1}C_{n/2}^{n/4 + x} C_{n/2}^{n/4 -
   x}\Sigma(x).
\]
По аналогии с доказательством леммы~\ref{ext1} можно при $x=0$
доказать существование таких решений
$(u^{\vec{x}}_1(1),\ldots,u^{\vec{x}}_{4}(1)),$ $\ldots,$
$(u^{\vec{x}}_1(m-2),\ldots,u^{\vec{x}}_{2^{m-1}}(m-2))$
систем~(\ref{system}) с $s=1$, $\ldots$, $s=m-2$ соответственно,
что для некоторой константы $C>0$ имеют место неравенства
\[
   \left|\frac{w^{\vec{x}}_j(s)}{2}-u^{\vec{x}}_j(s)\right|\le C.
\]
Очевидно, соответствующий этому решению член суммы из правой части
равенства~(\ref{sum2}) ограничен снизу величиной
$2^{(m-2)n}/n^{\tilde c}$ для некоторой константы $\tilde c>0$.
Следовательно,
\begin{equation}
   \Sigma(0)\ge 2^{(m-2)n} \big/ n^{\tilde c}.
\label{sigma_zero}
\end{equation}

По формуле Стирлинга $C_{n/2}^{n/4 + x} C_{n/2}^{n/4 - x} \le 2^n
e^{-8 x^2/n} (1 + o(1))$ при $|x| \le n^{0.6}$, а следовательно,
$C_{n/2}^{n/4 + x} C_{n/2}^{n/4 - x} \le 2^n e^{-8 n^{0.2}} (1 +
o(1))$ при $|x|
> n^{0.6}$. Отсюда в силу~(\ref{sigma_zero})
\[
    \sum_{|x| > n^{0.6}} C_{n/2}^{n/4 + x} C_{n/2}^{n/4 - x} \Sigma(x)\le
    n2^n e^{-8 n^{0.2}}2^{(m-2)n}= o(N_1 \Sigma(0)).
\]
Поэтому
\begin{equation}
    M_F \sim N \sum_{x = \left\lceil -n^{0.6} \right\rceil}^{\left\lfloor
    n^{0.6} \right\rfloor} C_{n/2}^{n/4 + x} C_{n/2}^{n/4 - x}
    \Sigma(x)\sim N \sum_{x = \left\lceil -n^{0.6} \right\rceil}^{\left\lfloor
    n^{0.6} \right\rfloor} C_{n/2}^{n/4 + x} C_{n/2}^{n/4 - x}
    \Sigma(0)
\label{constrict_x}
\end{equation}
по лемме~\ref{ext1}. Добавим в граф $F$ ребро между вершинами
$z_1$ и $z_2$ и обозначим полученный граф $F^+$. Поскольку
$M_{F^+} =N N_1 \Sigma(0)$ и
\[
    \sum_{x = \left\lceil -n^{0.6} \right\rceil}^{\left\lfloor n^{0.6} \right\rfloor}
    \left(C_{n/2}^{n/4 + x}\right)^2 \sim N,
\]
то $M_{F^+} \sim \frac{N_1}{N} M_F$ для любого $F\in\mathcal{M}_{m,L-1}$,
что доказывает утверждение леммы.

\subsection{Доказательство леммы \ref{ext2}}
\label{ext2_proof}
Будем считать, что индуцированный на множестве корней $R$ подграф
графа $H$ представляет собой клику. Тогда в силу
леммы~\ref{copies_lem}
\begin{equation}
    M_H \sim N^{d + k} \left(\frac{N_1}{N}\right)^{C_d^2 +
    l},\quad
    M_{H|_R} \sim N^{d} \left(\frac{N_1}{N}\right)^{C_d^2}.
\label{number_of_copies}
\end{equation}
Рассмотрим множество $X$ векторов $\vec{x}=(x_1,\ldots,x_{2^d})$,
для каждого из которых найдется последовательность $d$ вершин
$\mathbf{v}^1, \ldots, \mathbf{v}^d \in V$, образующих клику в
$G$, с таким вектором $\vec{x}$, определенным в~(\ref{x}). Для
каждого $\vec{x}\in X$ обозначим $T_d(\vec{x})$ количество
последовательностей $d$ вершин $\mathbf{v}^1, \ldots, \mathbf{v}^d
\in V$ с таким вектором $\vec{x}$. Тогда
\[
    M_H = \sum_{\vec{x}\in X} T_d(\vec{x}) M_{(R,
    H)}^{\vec{x}},\quad
    M_{H|_R} = \sum_{\vec{x}\in X} T_d(\vec{x}).
\]
Аналогично доказательству асимптотического
равенства~(\ref{constrict_x}), можно показать, что при
суммировании в двух последних равенствах по всем
$\vec{x}=(x_1,\ldots,$ $x_{2^d})$ из $X$ с условием $|x_i| \le
n^{0.6}$, $i \in \{1, \ldots, 2^d\}$, асимптотика
величины $M_H$ не поменяется. Поэтому в
силу~(\ref{number_of_copies})
\begin{multline*}
    M_H \sim \sum_{\substack{
    \vec{x} = \left(x_1, \ldots, x_{2^d}\right)\in X \\
    \forall i \in \left\{1, \ldots, 2^d\right\}\  |x_i| \le n^{0.6}
    }}
    T_d(\vec{x}) M_{(R, H)}^{\vec{x}} \sim \\
    \sum_{\substack{
    \vec{x} = \left(x_1, \ldots, x_{2^d}\right)\in X \\
    \forall i \in \left\{1, \ldots, 2^d\right\}\  |x_i| \le n^{0.6}
    }}
    T_d(\vec{x}) M_{(R, H)}' \sim
    N^d \left(\frac{N_1}{N}\right)^{C_d^2} M_{(R, H)}'.
\end{multline*}
Отсюда и из~(\ref{number_of_copies})
\[
    M_{(R, H)}' \sim N^k \left(\frac{N_1}{N}\right)^l,
\]
что и требовалось доказать.

\section{Доказательство теоремы \ref{copies_thrld}}
\label{copies_thrld_proof} Пусть $a$ --- количество автоморфизмов
графа $F$, а $\tilde M_F$
--- количество подграфов полного дистанционного графа $G$,
изоморфных графу $F$. Из леммы \ref{copies_lem} известно, что
\[
    M:=\tilde M_F \sim \frac1a N^{v(F)} \left(\frac{N_1}{N}\right)^{e(F)}.
\]
Пусть $F_1, \ldots, F_M$ --- это все копии $F$ в $G$.
Ясно, что $X_F$ можно представить в виде суммы
\[
    X_F = \sum_{i = 1}^M X_i,
\]
где $X_i$ --- индикатор того, что $F_i \subset G_p$.

Пользуясь линейностью математического ожидания, получим, что
\[
    \Expect X = \sum_{i=1}^M \Expect X_i \sim \frac1a N^{v(F)} \left(\frac{N_1}{N}\right)^{e(F)} p^{e(F)}.
\]

Пусть $H_0 \subseteq F$ --- подграф, на котором достигается
максимум плотности $\rho(H_0)$. Пусть $\tilde a$ --- количество
автоморфизмов графа $H_0$. Тогда
$\rho^{\max}(F)=\rho(H_0)=e(H_0)/v(H_0)$. Рассуждая аналогично,
найдем асимптотику математического ожидания $X_{H_0}$ (количества
копий $H_0$ в $G_p$):

\[
    \Expect X_{H_0}\sim \frac{1}{\tilde a} N^{v(H_0)} \left(\frac{N_1}{N}\right)^{e(H_0)} p^{e(H_0)}.
\]

Пусть
\[
    p \ll N^{-1/\rho(H_0)} \frac{N}{N_1}.
\]

Тогда $\Expect X_{H_0} \to 0$ и
\[
   \Prob(X > 0) \le \Prob(X_{H_0} >
   0)=\Prob(X_{H_0} \ge 1)\le\Expect X_{H_0} \to 0
\]
при $n \to \infty$ в силу неравенства Маркова.

Пусть теперь
\[
    p = \omega N^{-1/\rho(H_0)} \frac{N}{N_1},
\]
где $\omega = \omega(n) \to \infty$ при $n \to \infty$. Тогда для
любого графа $H\subseteq F$ имеем
\begin{equation}
   \Expect X_H\asymp
   N^{v(H)}N^{-\frac{e(H)}{\rho(H_0)}}\omega^{e(H)}\ge
   N^{v(H)}N^{-\frac{e(H)}{\rho(H)}}\omega^{e(H)}=\omega^{e(H)} \to \infty
\label{subgraph}
\end{equation}
при $n\rightarrow\infty$. Напомним, что
\[
    \Variance X_F = \sum_{ i = 1 }^M \Variance X_i + 2 \sum_{ 1 \le i < j \le M } \cov (X_i, X_j).
\]

Заметим, что случайные величины $X_i,X_j$ являются зависимыми в
том и только том случае, когда пересечение образов графа $F$ при
мономорфизмах $F\rightarrow F_i$, $F\rightarrow F_j$ содержит хотя
бы одно ребро. Кроме того, для любого подграфа $H\subset F$ число таких пар
подграфов $\{F_i, F_j\}$, что некоторый остовный подграф $F_i\cap
F_j$ изоморфен $H$, составляет
\[
    O\left(N^{2 v(F) - v(H)} \left(\frac{N_1}{N}\right)^{2 e(F) - e(H)}\right).
\]
Оценим дисперсию $X_F$:
\begin{multline*}
    \Variance X_F \le \Expect X_F +
    O\left(\sum_{\substack{H \subset F \\ e(H) \ge 1}} N^{2 v(F) - v(H)}
    \left(\frac{N_1}{N}\right)^{2 e(F) - e(H)} p^{2 e(F) - e(H)}\right) \le \\
    \Expect X_F + (\Expect X_F)^2 O\left(\sum_{\substack{H \subset F \\ e(H) \ge 1}}
    \frac{1}{N^{v(H)} \left(\frac{N_1}{N}\right)^{e(H)} p^{e(H)}}\right) = \\
    \Expect X_F + (\Expect X_F)^2 O\left(\sum_{\substack{H \subset F \\
    e(H) \ge 1}}
    \frac{1}{\Expect X_H}\right)=
    o\left((\Expect X_F)^2\right)
\end{multline*}
в силу~(\ref{subgraph}). Отсюда, воспользовавшись неравенством
Чебышева, получим, что
\begin{multline*}
    \Prob (X_F = 0) = \Prob (-X_F \ge 0) \le
    \Prob (|\Expect X_F - X_F| \ge \Expect X_F) \le \\
    \frac{ \Variance X_F }{ (\Expect X_F)^2 } \to 0 \text{ при $n \to \infty$.}
\end{multline*}
Применим это же неравенство для доказательства закона больших
чисел при $p \gg N^{-1/\rho^{\max}(F)} \frac{N}{N_1}$. Для любого
$\varepsilon
> 0$
\[
    \Prob\left(\left|\frac{X_F}{\Expect X_F} - 1 \right| < \varepsilon \right) =
    1 - \Prob(|X_F - \Expect X_F| \ge \varepsilon \Expect X_F) \ge
    1 - \frac{\Variance X_F}{\varepsilon^2 (\Expect X_F)^2} \to 1.
\]
Таким образом, теорема полностью доказана.

\section{Доказательство теоремы \ref{poisson}}
\label{poisson_proof} В силу леммы \ref{copies_lem}  количество
подграфов полного дистанционного графа $G$, изоморфных графу $F$,
равно
\[
    M \sim \frac1a N^k \left(\frac{N_1}{N}\right)^l.
\]
Отсюда в обозначениях из раздела~\ref{copies_thrld_proof}
\[
    \Expect X_F = \sum_{i=1}^M \Expect X_i = M p^l \sim c^l / a = \lambda.
\]

В силу теоремы 1.20 из \cite{bollobas} достаточно доказать
сходимость факториальных моментов
\[
    \Expect X^{\underline{j}} = \Expect X (X - 1) \ldots (X - j + 1)
\]
к $\lambda^j$ при $n \to \infty$ для всех $j \in \mathbb{N}$. Из
этого будет следовать утверждение теоремы.

Очевидно, $\Expect X^{\underline{j}}$ есть математическое ожидание
числа упорядоченных наборов $j$ различных копий графа $F$.
Покажем, что эта величина асимптотически совпадает с
математическим ожиданием $Y_j$
--- количества упорядоченных наборов $j$ копий $F$, множества вершин
которых не пересекаются. Заметим, что
\begin{multline}\label{mexpr}
    M \left(M - k N_1^{k-1}\right) \ldots \left(M - (j - 1) k N_1^{k-1}\right) p^{j l} \le \Expect Y_j \le \\
    M (M - 1) \ldots (M - j + 1) p^{j l},
\end{multline}
откуда следует, что
\[
    \Expect Y_j \sim M^j p^{j l} \sim \lambda^j.
\]

Остается доказать, что $\Expect Z_j \to 0$ при $n \to \infty$, где
$Z_j = X^{\underline{j}} - Y_j$.

Заметим сначала, что для двух таких различных графов $A$ и $B$ с $V(A) \cap V(B) \ne \varnothing$,
что $B$ изоморфен $F$ и в точности $t\ge 1$ вершин графа $B$ не
принадлежат $V(A)$, в силу строгой сбалансированности $F$ выполнено
$e\left(B|_{V(A \cap B)}\right) < \frac{(k - t) l}{k}$, из чего следует, что
\begin{equation}\label{lower1}
    \begin{aligned}
        e(A \cup B)=e(A)+e(B)-e\left(B|_{V(A \cap B)}\right)\ge e(A)+l-\frac{(k-t)l}{k}+\frac{1}{k}=\\
        e(A) + \frac{t l}{k} + \frac1k \text{ при $t<k$},
    \end{aligned}
\end{equation}
\begin{equation}\label{lower2}
    \begin{aligned}
        e(A \cup B)=e(A)+e(B)= e(A)+l\text{ при $t=k$}.
    \end{aligned}
\end{equation}
Пусть $F_{i_1}, \ldots, F_{i_j}$ --- подграфы $G$, изоморфные $F$,
с
\[
   \left|\bigcup_{v=1}^j V(F_{i_v})\right| = t < j k.
\]
Пусть для каждого $v\in\{2,\ldots,j\}$ ровно $t_v$ вершин подграфа
$F_{i_v}$ не принадлежат $\bigcup_{s=1}^{v-1} V(F_{i_s})$. Тогда,
очевидно, $\sum_{v=2}^j t_v = t - k$. Пусть сначала $t>k$. Без ограничения общности можно считать, что $0 < t_2 < k$, поэтому из~(\ref{lower1})
следует, что
\[
    e(F_{i_1} \cup F_{i_2}) \ge l + \frac{t_2 l}{k} + \frac1k.
\]
Поскольку в силу~(\ref{lower1})~и~(\ref{lower2})
\[
    e\left(\bigcup_{s=1}^v F_{i_s}\right)\ge e\left(\bigcup_{s=1}^{v-1} F_{i_s}\right) + \frac{t_v l}{k},
\]
получаем
\[
    e\left(\bigcup_{v=1}^{j} F_{i_v}\right) \ge \frac{t l}{k} + \frac1k.
\]
Заметим, что при $t=k$ для $e\left(\bigcup_{v=1}^{j} F_{i_v}\right)$ верна та же оценка.

Оценим математическое ожидание $Z_j$:
\[
    \Expect Z_j \le \sum_{t=k}^{k j - 1}
     O\left(N^t \left(\frac{N_1}{N}\right)^{\frac{t l}{k} + \frac1k} p^{\frac{t l}{k} + \frac1k}\right) =
     \sum_{t=k}^{k j - 1} O\left(\frac{N^t}{N^{\frac{k}{l} \left(\frac{t l}{k} + \frac1k\right)}}\right) =
     o(1).
\]

Таким образом, теорема доказана.

\section{Доказательство теоремы \ref{extdist}}
\label{extdist_proof} Будем для удобства обозначать
$\vec{\mathbf{v}} = (\mathbf{v}^1, \ldots, \mathbf{v}^d)$,
$\vec{\mathbf{w}} = (\mathbf{w}^1, \ldots, \mathbf{w}^k)$. Пусть
для любых $(\vec{\mathbf{v}}, \vec{\mathbf{w}}) \in \tilde
V_f^d \times V^k$
\begin{multline*}
    A_{\vec{\mathbf{v}} \vec{\mathbf{w}}}=\big\{G'\in\Omega_n^{\dist}:\,\,
    \text{ отображение $\lambda \colon (V(H),E(H)\setminus E(H|_R)) \to G'$,}\\
    \text{ такое что }
    \lambda(z_i) = \mathbf{v}^i, \, i\in\{1,\ldots,d\}, \, \lambda(y_j) = \mathbf{w}^j, \, j\in\{1,\ldots,k\},\\
    \text{ является мономорфизмом}\big\}.
\end{multline*}
Тогда, как легко видеть,
\begin{equation}
    \Ext_f^{\dist}(R, H) = \bigcap_{\vec{\mathbf{v}} \in \tilde V_f^d}
    \bigcup_{\vec{\mathbf{w}} \in V^k} A_{\vec{\mathbf{v}} \vec{\mathbf{w}}}.
\label{from_ext_to_a}
\end{equation}

Очевидно, что $\Prob(A_{\vec{\mathbf{v}} \vec{\mathbf{w}}}) =
p^l$. Более того, по лемме \ref{ext2} для любого
$\vec{\mathbf{v}}$ число различных $\vec{\mathbf{w}}$, для которых
$G\in A_{\vec{\mathbf{v}} \vec{\mathbf{w}}}$, асимптотически равно
$M_{(R, H)}' := N^k \left(\frac{N_1}{N}\right)^l$.

Будем называть последовательности
$\vec{\mathbf{w}}=(\mathbf{w}^1,\ldots,\mathbf{w}^d)$ и
$\vec{\mathbf{w}}'=((\mathbf{w}')^1,\ldots,$ $(\mathbf{w}')^d)$
эквивалентными, если они совпадают как множества и существует
автоморфизм $H$, оставляющий на месте корни, который переводит
каждое $y_i$ в $y_{\sigma(i)}$, $i\in\{1,\ldots,k\}$, где
перестановка $\sigma$ на $\{1,\ldots,d\}$ определена следующим
образом: $(\mathbf{w}')^i=\mathbf{w}^{\sigma(i)}$ для всех
$i\in\{1,\ldots,d\}$. Таким образом, последовательности
$\vec{\mathbf{w}}$ разбиваются на $c_1$ классов эквивалентности.
Поскольку для эквивалентных $\vec{\mathbf{w}}$ и
$\vec{\mathbf{w}}'$ события $A_{\vec{\mathbf{v}}
\vec{\mathbf{w}}}$ и $A_{\vec{\mathbf{v}} \vec{\mathbf{w}}'}$
совпадают, будем рассматривать для каждого $\vec{\mathbf{v}}$ по
одному представителю из каждого класса эквивалентности (множество
которых обозначим $V(\vec{\mathbf{v}})$).

Пусть
\[
    B_{\vec{\mathbf{v}}} = \overline{\bigcup_{\vec{\mathbf{w}}\in V(\vec{\mathbf{v}})}
    A_{\vec{\mathbf{v}} \vec{\mathbf{w}}}} = \bigcap_{\vec{\mathbf{w}}\in V(\vec{\mathbf{v}})}
    \overline{A_{\vec{\mathbf{v}} \vec{\mathbf{w}}}}.
\]

Пользуясь корреляционным неравенством из \cite{spencer}, получаем
\begin{equation}
    \prod_{\vec{\mathbf{w}}\in V(\vec{\mathbf{v}})}
    \Prob\left(\overline{A_{\vec{\mathbf{v}} \vec{\mathbf{w}}}}\right)
    \le \Prob(B_{\vec{\mathbf{v}}}) \le
    \left(\prod_{\vec{\mathbf{w}}\in V(\vec{\mathbf{v}})}
    \Prob\left(\overline{A_{\vec{\mathbf{v}} \vec{\mathbf{w}}}}\right)\right)
    \exp\left(2 \sum \Prob(A_{\vec{\mathbf{v}} \vec{\mathbf{w}}} \cap
    A_{\vec{\mathbf{v}} \vec{\mathbf{w}}'})\right),
\label{b_estimation}
\end{equation}
где суммирование под экспонентой ведется по всем $\vec{\mathbf{w}}
= (\mathbf{w}^1, \ldots, \mathbf{w}^k),$ $\vec{\mathbf{w}}' =
((\mathbf{w}')^1, \ldots, (\mathbf{w}')^k)$ с $\{\mathbf{w}^1,
\ldots, \mathbf{w}^k\} \cap \{(\mathbf{w}')^1, \ldots,
(\mathbf{w}')^k\} \ne \varnothing$.

Пусть
\begin{equation}
   p^l=cc_1d\ln N N^{-k}\left(\frac{N}{N_1}\right)^l.
\label{p_def}
\end{equation}
Тогда
\begin{multline}
    \prod_{\vec{\mathbf{w}}\in V(\vec{\mathbf{v}})}
    \Prob\left(\overline{A_{\vec{\mathbf{v}} \vec{\mathbf{w}}}}\right) = (1 - p^l)^
    {\frac{1}{c_1} N^k \left(\frac{N_1}{N}\right)^l (1 + h_{\vec{\mathbf{v}}}(n))} =
    e^{- \frac{1}{c_1} p^l N^k \left(\frac{N_1}{N}\right)^l (1 + g_{\vec{\mathbf{v}}}(n))} = \\
    e^{- c d \ln N (1 + g_{\vec{\mathbf{v}}}(n))} = N^{- c d (1 + g_{\vec{\mathbf{v}}}(n))},
\label{a_estimation}
\end{multline}
где $h_{\vec{\mathbf{v}}}(n) = O\left(f(n) n^{-0.8}\right)$, $g_{\vec{\mathbf{v}}}(n) = O\left(f(n) n^{-0.8}\right)$ при $n\rightarrow\infty$.

Разобьем все такие пары $(\vec{\mathbf{w}}, \vec{\mathbf{w}}')$ в
суммировании в~(\ref{b_estimation}), что $G\in
A_{\vec{\mathbf{v}}\vec{\mathbf{w}}}$, на классы $W_t$,
$t\in\{1,\ldots,k\}$, в которых
\[
   \left|\{\mathbf{w}^1, \ldots, \mathbf{w}^k\} \cap \{(\mathbf{w}')^1,
   \ldots, (\mathbf{w}')^k\}\right| = t.
\]
Из лемм~\ref{ext1}--\ref{ext2} следует, что
\[
   |W_t|=O\left(N^{2 k - t} \left(\frac{N_1}{N}\right)^{2 l - e_t}\right)
\]
равномерно по $\vec{\mathbf{v}} \in \tilde V_f^d$, где $e_t$ --- наибольшее по всем парам $(\vec{\mathbf{w}},
\vec{\mathbf{w}}')\in W_t$ количество ребер из
$E(G|_{\mathbf{v}^1,\ldots,\mathbf{v}^d,\mathbf{w}^1,\ldots,\mathbf{w}^k})\cap
E(G|_{\mathbf{v}^1,\ldots,\mathbf{v}^d,(\mathbf{w}')^1,\ldots,(\mathbf{w}')^k})
\setminus E(G|_{\mathbf{v}^1,\ldots,\mathbf{v}^d})$. Поскольку
сеть $(R, H)$ строго сбалансирована, $e_t < t l / k$, если $t <
k$. Если же $t = k$, то $e_t < l$, так как иначе расширения
$\vec{\mathbf{w}}$ и $\vec{\mathbf{w}}'$ были бы эквивалентны.
Поэтому в силу~(\ref{p_def})
\begin{multline*}
    \sum_{(\vec{\mathbf{w}},\vec{\mathbf{w}}')\in W_t}
    \Prob(A_{\vec{\mathbf{v}} \vec{\mathbf{w}}} \cap
    A_{\vec{\mathbf{v}} \vec{\mathbf{w}}'}) =
    O\left(N^{2 k - t} \left(\frac{N_1}{N}\right)^{2 l - e_t}
    p^{2 l - e_t}\right) = \\
    O\left(\frac{\ln^2 N}{N^t
    \left(\frac{N_1}{N}\right)^{e_t} p^{e_t}}\right) = o(1)
\end{multline*}
равномерно по $\vec{\mathbf{v}} \in \tilde V_f^d$. Отсюда $\sum_{t=1}^k\sum_{(\vec{\mathbf{w}},\vec{\mathbf{w}}')\in
W_t} \Prob(A_{\vec{\mathbf{v}} \vec{\mathbf{w}}} \cap
A_{\vec{\mathbf{v}} \vec{\mathbf{w}}'}) = o(1)$ равномерно по $\vec{\mathbf{v}} \in \tilde V_f^d$, и
\begin{equation}
    \Prob(B_{\vec{\mathbf{v}}}) = N^{- c d (1 + g_{\vec{\mathbf{v}}}(n))}(1+o(1))
\label{prob_b}
\end{equation}
равномерно по $\vec{\mathbf{v}} \in \tilde V_f^d$, так как справедливы
оценки~(\ref{b_estimation})~и~(\ref{a_estimation}).

Теперь легко видеть, что при $c > 1$ в силу~(\ref{from_ext_to_a})
\[
    \Prob\left(\overline{\Ext_f^{\dist}(R, H)}\right) \le N^d N^{- c d (1 + o(1))} \to 0
    \text{ при }n\rightarrow\infty.
\]

Рассмотрим теперь случай $c < 1$. Пусть случайная величина
$X_{\vec{\mathbf{v}}}$ есть индикатор события
$B_{\vec{\mathbf{v}}}$. Будем называть $\vec{\mathbf{v}} =
(\mathbf{v}^1, \ldots, \mathbf{v}^d)$ и $\vec{\mathbf{v}}' =
((\mathbf{v}')^1, \ldots,$ $(\mathbf{v}')^d)$ эквивалентными, если
они совпадают как множества и существует автоморфизм $H$, который
переводит каждое $z_i$ в $z_{\sigma(i)}$, $i\in\{1,\ldots,d \}$,
где перестановка $\sigma$ на $\{1,\ldots,d\}$ определена следующим
образом: $(\mathbf{v}')^i=\mathbf{v}^{\sigma(i)}$ для всех
$i\in\{1,\ldots,d\}$. Очевидно, для эквивалентных
$\vec{\mathbf{v}}, \vec{\mathbf{v}}'$ выполнено
$X_{\vec{\mathbf{v}}} = X_{\vec{\mathbf{v}}'}$. Как и ранее, будем
рассматривать только по одному представителю из каждого класса
эквивалентности (обозначим множество таких представителей $U$).
Пусть $X = \sum_{\vec{\mathbf{v}}\in U} X_{\vec{\mathbf{v}}}$.
Заметим, что $X = 0$ тогда и только тогда, когда выполнено
свойство $\Ext_f^{\dist}(R, H)$. В силу~(\ref{prob_b}) имеем
\[
    \Expect X = \sum_{\vec{\mathbf{v}}\in U} \Expect X_{\vec{\mathbf{v}}}
    \ge N^d N^{- c d (1 + o(1))} \to \infty\text{ при }n\rightarrow\infty.
\]
Очевидно,
\[
    \Variance X = \sum_{\vec{\mathbf{v}}\in U} \Variance X_{\vec{\mathbf{v}}} +
    \sum_{\vec{\mathbf{v}} \ne \vec{\mathbf{v}}'}
    \cov(X_{\vec{\mathbf{v}}}, X_{\vec{\mathbf{v}}'}).
\]
Оценим слагаемое $\Expect(X_{\vec{\mathbf{v}}}
X_{\vec{\mathbf{v}}'})$ в выражении для ковариации, пользуясь
корреляционным неравенством из \cite{spencer}. Так как при
некотором выборе множеств
$V(\vec{\mathbf{v}}),V(\vec{\mathbf{v}}')$ множества таких
$\vec{\mathbf{w}}\in V(\vec{\mathbf{v}})$ и $\vec{\mathbf{w}}\in
V(\vec{\mathbf{v'}})$, что
$\mathbf{w}^i\notin\{\mathbf{v}^j,(\mathbf{v}')^j\}$ для всех
$i\in\{1,\ldots,k\}$ и $j\in\{1,\ldots,d\}$, совпадают, то
\begin{multline*}
    \Expect(X_{\vec{\mathbf{v}}} X_{\vec{\mathbf{v}}'}) =
    \Prob(B_{\vec{\mathbf{v}}} \cap B_{\vec{\mathbf{v}}'}) =
    \Prob\left(\left(\bigcap_{\vec{\mathbf{w}}\in V(\vec{\mathbf{v}})}
    \overline{A_{\vec{\mathbf{v}} \vec{\mathbf{w}}}}\right) \cap
    \left(\bigcap_{\vec{\mathbf{w}}\in V(\vec{\mathbf{v}}')}
    \overline{A_{\vec{\mathbf{v}}' \vec{\mathbf{w}}}}\right)\right) \le \\
    \Prob\left(\left(\bigcap_{\vec{\mathbf{w}}\in V(\vec{\mathbf{v}}) \colon
    \forall i\in\{1,\ldots,k\}\forall j\in\{1,\ldots,d\}\,
    \mathbf{w}^i\notin\{\mathbf{v}^j,(\mathbf{v}')^j\}}
    \overline{A_{\vec{\mathbf{v}} \vec{\mathbf{w}}}}\right) \cap \right. \\
    \left. \left(\bigcap_{\vec{\mathbf{w}}\in V(\vec{\mathbf{v}}') \colon
    \forall i\in\{1,\ldots,k\}\forall j\in\{1,\ldots,d\}\,
    \mathbf{w}^i\notin\{\mathbf{v}^j,(\mathbf{v}')^j\}}
    \overline{A_{\vec{\mathbf{v}}' \vec{\mathbf{w}}}}\right)\right) \le \\
    \prod_{\vec{\mathbf{w}}\in V(\vec{\mathbf{v}}), \vec{\mathbf{w}}'\in V(\vec{\mathbf{v}}') \colon
    \forall i\in\{1,\ldots,k\}\forall j\in\{1,\ldots,d\}\,
    \{\mathbf{w}^i, (\mathbf{w}')^i\} \cap \{\mathbf{v}^j,(\mathbf{v}')^j\} = \varnothing}
    \Prob\left(\overline{A_{\vec{\mathbf{v}} \vec{\mathbf{w}}}}\right)
    \Prob\left(\overline{A_{\vec{\mathbf{v}}'
    \vec{\mathbf{w}}'}}\right)\times\\
    \exp\left(2 \left(\sum\left( \Prob(A_{\vec{\mathbf{v}} \vec{\mathbf{w}}} \cap
    A_{\vec{\mathbf{v}} \vec{\mathbf{w}}'})+
    \Prob(A_{\vec{\mathbf{v}}' \vec{\mathbf{w}}} \cap
    A_{\vec{\mathbf{v}}' \vec{\mathbf{w}}'})\right)+
    \sum \Prob(A_{\vec{\mathbf{v}} \vec{\mathbf{w}}} \cap
    A_{\vec{\mathbf{v}}' \vec{\mathbf{w}}'})\right)\right),
\end{multline*}
где суммирование в первой сумме ведется по наборам $\vec{\mathbf{w}} =
(\mathbf{w}^1, \ldots, \mathbf{w}^k),$ $\vec{\mathbf{w}}' =
((\mathbf{w}')^1, \ldots, (\mathbf{w}')^k)$ с $\{\mathbf{w}^1,
\ldots, \mathbf{w}^k\} \cap \{(\mathbf{w}')^1, \ldots,
(\mathbf{w}')^k\} \ne \varnothing$, а во второй --- по
$\vec{\mathbf{w}},$ $\vec{\mathbf{w}}'$ с пересекающимися
множествами ребер
$E(G|_{\{\mathbf{v}^1,\ldots,\mathbf{v}^d,\mathbf{w}^1,\ldots,\mathbf{w}^k\}})\setminus
E(G|_{\{\mathbf{v}^1,\ldots,\mathbf{v}^d\}})$ и
$E(G|_{\{(\mathbf{v}')^1,\ldots,(\mathbf{v}')^d,(\mathbf{w}')^1,\ldots,(\mathbf{w}')^k\}})\setminus
E(G|_{\{(\mathbf{v}')^1,\ldots,(\mathbf{v}')^d\}})$ (причем в
обеих суммах $\vec{\mathbf{w}}$ и $\vec{\mathbf{w}}'$ удовлетворяют
условию
$\{\mathbf{w}^i,(\mathbf{w}')^i\}\cap\{\mathbf{v}^j,(\mathbf{v}')^j\}=\varnothing$
для всех $i\in\{1,\ldots,k\}$ и $j\in\{1,\ldots,d\}$). Выше мы уже
доказали, что
\[
    \sum \Prob(A_{\vec{\mathbf{v}} \vec{\mathbf{w}}} \cap A_{\vec{\mathbf{v}} \vec{\mathbf{w}}'}) = o(1),\quad
    \sum \Prob(A_{\vec{\mathbf{v}}' \vec{\mathbf{w}}} \cap A_{\vec{\mathbf{v}}' \vec{\mathbf{w}}'}) =
    o(1)
\]
равномерно по $\vec{\mathbf{v}}, \vec{\mathbf{v}}' \in \tilde V_f^d$. Таким образом, нужно показать, что также равномерно по $\vec{\mathbf{v}}, \vec{\mathbf{v}}' \in \tilde V_f^d$ выполнено
\begin{equation}
    \sum \Prob(A_{\vec{\mathbf{v}} \vec{\mathbf{w}}} \cap
    A_{\vec{\mathbf{v}}' \vec{\mathbf{w}}'}) = o(1).
\label{final_sum}
\end{equation}
Опять разобьем все такие пары $(\vec{\mathbf{w}},
\vec{\mathbf{w}}')$ в суммировании, что $G\in
A_{\vec{\mathbf{v}}\vec{\mathbf{w}}}$, на классы $W_t$,
$t\in\{1,\ldots,k\}$, в которых
\[
   \left|\{\mathbf{w}^1, \ldots, \mathbf{w}^k\} \cap \{(\mathbf{w}')^1,
   \ldots, (\mathbf{w}')^k\}\right| = t.
\]
В силу наших ограничений на пределы суммирования
в~(\ref{final_sum}) для любой пары $(\vec{\mathbf{w}},
\vec{\mathbf{w}}')\in W_t$ выполнено
\[
   \left(\{\mathbf{v}^1,\ldots,\mathbf{v}^d\} \cup
    \{(\mathbf{v}')^1,\ldots,(\mathbf{v}')^d\}\right) \cap
    \left(\{\mathbf{w}^1,\ldots,\mathbf{w}^k\} \cup
    \{(\mathbf{w}')^1,\ldots,(\mathbf{w}')^k\}\right) =\varnothing.
\]
Поэтому из лемм~\ref{ext1}--\ref{ext2} следует
\[
   |W_t|=O\left(N^{2 k - t} \left(\frac{N_1}{N}\right)^{2 l -
   e_t}\right),
\]
где $e_t < t l / k$, равномерно по $\vec{\mathbf{v}}, \vec{\mathbf{v}}' \in \tilde V_f^d$. Таким образом, снова получаем, что для
каждого $t\in\{1,\ldots,k\}$
\[
   \sum_{(\vec{\mathbf{w}},
   \vec{\mathbf{w}}')\in W_t} \Prob(A_{\vec{\mathbf{v}}
   \vec{\mathbf{w}}} \cap
   A_{\vec{\mathbf{v}}' \vec{\mathbf{w}}'})=o(1)
\]
равномерно по $\vec{\mathbf{v}}, \vec{\mathbf{v}}' \in \tilde V_f^d$, поэтому по аналогии с~(\ref{a_estimation})
\begin{multline*}
    \Expect(X_{\vec{\mathbf{v}}} X_{\vec{\mathbf{v}}'}) \le \\
    \prod_{\vec{\mathbf{w}}\in V(\vec{\mathbf{v}}), \vec{\mathbf{w}}'\in V(\vec{\mathbf{v}}') \colon
    \forall i\in\{1,\ldots,k\}\forall j\in\{1,\ldots,d\}\,
    \{\mathbf{w}^i, (\mathbf{w}')^i\} \cap \{\mathbf{v}^j,(\mathbf{v}')^j\} = \varnothing}
    \Prob\left(\overline{A_{\vec{\mathbf{v}} \vec{\mathbf{w}}}}\right) \times \\
    \Prob\left(\overline{A_{\vec{\mathbf{v}}'
    \vec{\mathbf{w}}'}}\right)(1+o(1))\le
    \left(1-p^l\right)^{\frac{1}{c_1} N^k \left(\frac{N_1}{N}\right)^l (1+h_{\vec{\mathbf{v}}}(n)) - k d N^{k-1}} \times \\ \left(1-p^l\right)^{\frac{1}{c_1} N^k \left(\frac{N_1}{N}\right)^l (1+h_{\vec{\mathbf{v}}'}(n)) - k d N^{k-1}} (1+o(1)) \le \\
    N^{- c d (2 + g_{\vec{\mathbf{v}}}(n) + g_{\vec{\mathbf{v}}'}(n))} (1 + o(1))
\end{multline*}
равномерно по $\vec{\mathbf{v}}, \vec{\mathbf{v}}' \in \tilde V_f^d$. Отсюда в силу~(\ref{prob_b})
\begin{multline*}
    \cov(X_{\vec{\mathbf{v}}}, X_{\vec{\mathbf{v}}'}) =
    \Expect(X_{\vec{\mathbf{v}}} X_{\vec{\mathbf{v}}'}) -
    \Expect X_{\vec{\mathbf{v}}} \Expect X_{\vec{\mathbf{v}}'} \le
    N^{- c d (2 + g_{\vec{\mathbf{v}}}(n) + g_{\vec{\mathbf{v}}'}(n))} (1 + o(1)) - \\
    N^{- c d (2 + g_{\vec{\mathbf{v}}}(n) + g_{\vec{\mathbf{v}}'}(n))} =
    o\left(N^{- c d (2 + g_{\vec{\mathbf{v}}}(n) + g_{\vec{\mathbf{v}}'}(n))}\right)=o(\Expect X_{\vec{\mathbf{v}}} \Expect X_{\vec{\mathbf{v}}'})
\end{multline*}
равномерно по $\vec{\mathbf{v}}, \vec{\mathbf{v}}' \in \tilde V_f^d$. Поэтому
\[
    \Variance X \le \Expect X + \sum_{\vec{\mathbf{v}} \ne \vec{\mathbf{v}}'} o(\Expect X_{\vec{\mathbf{v}}} \Expect X_{\vec{\mathbf{v}}'}) =
    o\left((\Expect X)^2\right).
\]
Пользуясь неравенством Чебышева, получаем, что
\[
    \Prob (X = 0) = \Prob (-X \ge 0) \le \Prob (|\Expect X - X| \ge \Expect X) \le \frac{ \Variance X }{ (\Expect X)^2 } \to 0.
\]

Теорема доказана.

\end{document}